\documentclass[10pt,notitlepage]{amsart}
\usepackage{bezier,graphicx} 
\usepackage{amscd}
\usepackage{enumerate}

\def\cocoa{{\hbox{\rm C\kern-.13em o\kern-.07em C\kern-.13em o\kern-.15em A}}}

\def\TS{\mathcal{X}}

\def\w{wt(2)}

\def\ox{\otimes}
\def\wh{\widehat}

\advance\hoffset by 0.4 truein
\def\lra{\longrightarrow}

\renewcommand{\:}{\colon}
\newcommand{\QQ}{{\mathbb Q}}
\newcommand{\CC}{{\mathbb C}}
\newcommand{\ol}{\overline}
\newcommand{\wt}{\widetilde}
\newcommand{\cal}{\mathcal}
 \font\tenbi=cmmi14
 \font\sevenbi=cmmi10 \font\fivebi=cmmi7
 \newfam\bifam \def\bi{\fam\bifam} \textfont\bifam=\tenbi
 \scriptfont\bifam=\sevenbi \scriptscriptfont\bifam=\fivebi
 \mathchardef\variablemega="7121 \def\w{{\bi\variablemega}}

\newtheorem{theorem}[subsection]{Theorem}
\newtheorem{proposition}[subsection]{Proposition}
\newtheorem{lemma}[subsection]{Lemma}
\theoremstyle{definition}  
\newtheorem{claim}[subsection]{}
\newtheorem{remark}[subsection]{Remark}

\font\bigbf=cmbx16

\font\largebf=cmbx12

\begin{document}

\

\title[{\bf Limits of special Weierstrass points}\hfill]
{\baselineskip=12pt \bigbf Limits of special 
Weierstrass \\ \  \\ points}

\author[{\bf\hfill Cumino, Esteves and Gatto}]
{\vskip-1.3cm}

\maketitle

\

\begin{center}
{\largebf C. Cumino$^{\text{\bf a}}$, 
E. Esteves$^{\text{\bf b}}$ and 
L. Gatto$^{\text{\bf a}}$
\footnote{Work partially sponsored by MURST 
(Progetto Nazionale ``Geometria delle Variet\`a Proiettive'' Coordinatore 
Sandro Verra), and supported by 
GNSAGA-INDAM. The second author was also supported
by CNPq, Proc. 478625/03-0 and 301117/04-7, 
and CNPq/FAPERJ, Proc. E-26/171.174/2003.}}
\end{center}

\vskip0.4cm

\begin{center}
$^{\text{\rm a}}$ 
Dipartimento di Matematica, Politecnico di Torino,\\
C.so Duca degli Abruzzi 24, 10129 Torino, Italy\\
\  \\
$^{\text{\rm b}}$ 
Instituto Nacional de Matem\'atica Pura e Aplicada,\\ 
Estrada Dona Castorina 110, 22460-320 Rio de Janeiro, 
Brazil
\end{center}

\vskip0.6cm

\section{Introduction}

\begin{claim} (\emph{Our main result}) 
Let $X\cup_PY$ be the union of two general connected, 
smooth, nonrational curves $X$ and $Y$ intersecting 
transversally at a point $P$. 
Assume that $P$ is a general point of $X$ or of $Y$.
Our main result is 
Theorem~\ref{mainthm}, which, in a simplified way, 
says:

\medskip

{\em Let $Q\in X$. Then $Q$ is the 
limit of special Weierstrass points on a 
family of smooth 
curves degenerating to $C$ if and only if $Q\neq P$ and 
either of the following conditions hold: $Q$ is 
a special ramification point of the linear system 
$|K_X+(g_Y+1)P|$, or $Q$ is a ramification point of the 
linear system $|K_X+(g_Y+1+j)P|$ for $j=\pm 1$ 
and $P$ is a Weierstrass point of $Y$.}
\medskip

Above, $g_Y$ stands for the genus of $Y$ and 
$K_X$ for a canonical divisor of $X$.

As an application, 
we use Theorem~\ref{mainthm} to recover in a 
unified way computations made by Diaz and Cukierman 
of certain divisor classes in the moduli spaces of 
stable curves; see Subsections 1.2 and 1.3 and 
Sections 5 and 6.
\end{claim}

\begin{claim} (\emph{Motivation}) 
In order to understand how the above result fits in 
the literature on the subject, we must recall that in 
the last two decades several papers on limits of 
Weierstrass points and linear series on stable curves 
appeared. The investigations about these topics 
were initially aimed to prove existence theorems (about, 
e.g., distinguished linear series on smooth curves) or 
to do enumerative geometry, 
in the sense of \cite{MumfEnum}, on the moduli space of 
genus-$g$ stable curves $\ol M_g$. 
For instance, in 
the beginning of the eighties, Harris and Mumford 
\cite{HaMu} proved 
that the moduli space $\overline{M}_g$ is of general 
type for $g$ odd and $g\geq 23$, doing computations on 
${\rm Pic}(\overline{\cal M}_g)$, 
the Picard group of the moduli functor. 

The same techniques were successfully used by Diaz 
\cite{Diazexc} to compute the class 
$\ol E_{g,-1}$ (there named $\overline{D_{g-1}}$) 
of the closure of the locus of curves having 
an exceptional (here called {\em special}) Weierstrass 
point of type $g-1$; see Subsection~\ref{srl} for a 
precise definition. A Weierstrass point $Q$ on a 
smooth curve of genus $g$ is said to be of type $g-1$ if 
$\dim|(g-1)Q|\geq 1$, and of type $g+1$ if 
$\dim|(g+1)Q|\geq 2$. Diaz computed the class 
$\ol E_{g,-1}$ by intersecting it with 
certain test curves entirely contained in the boundary 
of $M_g$ in $\overline{M}_g$. This way he got relations 
among the coefficients of the expression of 
$\ol E_{g,-1}$ in 
terms of generators of ${\rm Pic}(\overline{\cal M}_g)$. 
These test curves were induced by one-parameter 
families ${\cal F}_i\to X_i$ of curves given as follows: 
start with a general smooth curve $X_i$ of genus $g-i$, 
for each $i=1,\dots,g-1$, and a 
general smooth pointed curve $(Y_i,B_i)$ of genus $i$; 
then the fiber $({\cal F}_i)_P$ over $P\in X_i$ is 
$X_i\cup_PY_i$, the point $B_i\in Y_i$ 
being identified with $P\in X_i$. 
This can be seen as a curve in $\overline{M}_g$ via a 
nonconstant map $\gamma_i\:X_i\to\overline{M}_g$.

The crux of Diaz's method was to evaluate 
$\int_X\gamma_i^*\ol E_{g,-1}$, which amounts to 
knowing, 
with multiplicities, for how many pairs $(P,R)$ with 
$P\in X_i$ and $R\in X_i\cup_P Y_i$ there is a family of 
smooth curves 
degenerating to $X_i\cup_PY_i$ with Weierstrass points 
of type $g-1$ converging to $R$. This was done in 
\cite{Diazexc} by 
using the theory of admissible coverings 
introduced and developed in \cite{HaMu}. So half of 
our Theorem~\ref{mainthm} is in \cite{Diazexc}.

After Diaz's work, it was natural to ask what the 
limits of special Weierstrass points of type $g+1$ are, 
the other half of Theorem~\ref{mainthm}. 
In fact, soon afterwards, 
Cukierman \cite{cukie} computed the class 
$\ol E_{g,1}$ of the closure of the locus of 
curves having a Weierstrass point of type $g+1$; see 
Subsection~\ref{srl}. However, his method 
was not based on 
test curves, but on a Hurwitz formula 
with singularities. (He used Diaz's result as well.) 
Also, the theory of admissible coverings could not be 
effectively used, as the 
condition defining $\ol E_{g,1}$ is not about the 
existence of a pencil, but of a net. Of course, 
once we have 
an expression for $\ol E_{g,1}$ in terms of the 
generators of ${\rm Pic}(\overline{\cal M}_g)$, we can 
evaluate it along the $\gamma_i$. But we 
cannot infer what the limits of Weierstrass 
points of type $g+1$ on 
$X_i\cup_PY_i$ are just from 
their number. 

Our Theorem~\ref{mainthm} fills this 
gap.
To show the 
``only if'' 
part of it is not hard. 
To show the 
``if'' 
part we use limit 
linear series on two-parameter families of curves, 
instead of admissible coverings.
\end{claim}

\begin{claim} (\emph{Application})
Our Theorem~\ref{mainthm} can be used to 
compute the  
classes $\ol E_{g,-1}$ and 
$\ol E_{g,1}$ in a unified and conceptually simpler way. 
Also, 
we do not need to worry about multiplicities, an usual 
nuisance of the method of test curves.

In brief, here is how.
First of all, 
we consider another divisor class on $\ol M_g$, 
the class $\ol{SW}_g$ of the closure of the locus of 
curves having a special Weierstrass point, either of 
type $g-1$ or of type $g+1$; see Subsection \ref{srl} 
for a more precise definition. It turns out that 
$\ol{SW}_g$ is much easier to compute. An expression for 
it, in terms of the generators of 
${\rm Pic}(\overline{\cal M}_g)$, appeared already 
in \cite{mathscand}, though there are multiplicity 
issues due to the method of test curves. 

Here we compute $\ol{SW}_g$ in Theorem \ref{formSW} 
directly, by intersecting $\ol{SW}_g$ with a general 
curve in $\ol M_g$. No multiplicity issues arise. Of 
vital importance in this computation is 
Theorem~\ref{irredlim}, which essentially 
describes the limits of Weierstrass points on a 
general irreducible uninodal curve. This 
description is much finer 
than the one found in \cite{Diazexc}, Thm. A2.1, p. 60, 
for instance. For the proof of Theorem~\ref{irredlim} we 
use the theory of limit linear series for curves that 
are not of compact type, developed in \cite{EST2}. 

Then we show that $\ol{SW}_g=\ol E_{g,-1}+\ol E_{g,1}$. 
This follows from our Proposition \ref{schthunion}. 
This is something to be expected, from a purely 
set-theoretic point of view, but nevertheless, because 
of multiplicity issues, is not immediate and 
had to be proved. 

Now we use the test curves given by the $\gamma_i$. 
Having the expression for $\ol{SW}_g$ allows us to 
compute $\int_X\gamma_i^*\ol{SW}_g$, which gives us the 
sum
$$
\int_X\gamma_i^*\ol E_{g,-1}+
\int_X\gamma_i^*\ol E_{g,1}
$$
for each $i=1,\dots,g-1$. For each $j=-1,1$, 
let $e_{i,j}$ denote the number of pairs $(P,R)$ 
with $P\in X_i$ and $R\in X_i\cup_PY_i$ 
such that there is a 
family of smooth curves degenerating to 
$X_i\cup_PY_i$ with 
Weierstrass points of type $g+j$ converging to $R$. 
Theorem \ref{mainthm} tells us what these pairs are. 
Their number is computed in \cite{cumestgat1}, Thm. 5.6. 
Then
\begin{equation}\label{iE}
\int_X\gamma_i^*\ol E_{g,j}\geq e_{i,j}.
\end{equation}
In principle, the inequality may 
be strict
because of multiplicity issues. However, a 
simple computation yields
$$
\int_X\gamma_i^*\ol{SW}_g=e_{i,-1}+e_{i,1}.
$$
Thus equality holds in \eqref{iE}. From this equality, 
for each $j=-1,1$ and each $i=1,\dots,g-1$, the classes 
$\ol E_{g,-1}$ and $\ol E_{g,1}$ are computed in the 
usual way, like in \cite{Diazexc}; see 
Subsection \ref{testfamily} for more details. 
\end{claim}

\begin{claim} (\emph{Layout}) In Section 2 we present a 
few preliminaries on ramification schemes, deformations 
of curves and limit linear series. In Section 3 we 
describe limits of special Weierstrass points on 
reducible uninodal curves, and in Section 4 we describe 
limits of Weierstrass points on irreducible uninodal 
curves. In Section 5 we compute $\ol{SW}_g$ and in 
Section 6 we compute $\ol{E}_{g,-1}$ and 
$\ol{E}_{g,+1}$. 
\end{claim}
  
\begin{claim} (\emph{Acknowledgments}) The authors wish 
to thank Nivaldo Medeiros for discussions on related 
topics and acknowledge the use of \cocoa \cite{cocoa} 
for some of the computations.
\end{claim}

\section{Preliminaries}

\begin{claim}\label{prel1} 
({\em Ramification points}) A ({\em nodal}) 
{\em curve} is a connected, reduced, projective scheme of 
dimension 1 over $\mathbb C$ whose 
only singularities are nodes, 
i.e. ordinary double points. 
The canonical sheaf, or dualizing sheaf of a curve $C$ 
will be denoted $\w_C$. By the hypothesis on the 
singularities of $C$, the sheaf $\w_C$ is a line bundle. 
The (arithmetic) genus of $C$, i.e. $h^0(C,\w_C)$, will 
be denoted $g_C$. 

Let $C$ be a smooth curve, and $V$ a linear system of 
dimension $r$ of sections of a line bundle $L$ on $C$. 
For each $P\in C$ and each nonnegative integer $a$, let 
$V(-aP)\subseteq V$ be the linear subsystem of 
sections of $V$ vanishing at $P$ with multiplicity 
at least $a$. We call $P$ a 
{\em ramification point} of $V$ 
if $\dim V(-rP)\geq 1$; otherwise we call $P$ an 
{\em ordinary point} of $V$.  
A ramification point $P$ 
of $V$ is said to be {\em special of type $r-1$} if 
$\dim V(-(r-1)P)\geq 2$, and 
{\em special of type $r+1$} if 
$\dim V(-(r+1)P)\geq 1$. 
A {\em special ramification point} of $V$ is a 
special ramification point of type $r-1$ or $r+1$.

The orders of vanishing at $P$ of the 
sections of $L$ in $V$ 
can be ordered increasingly. We call 
this increasing sequence 
the {\em order sequence} of $V$ at $P$. 
The order sequence is 
$0,1,\dots,r-1$ if and only if $P$ is an ordinary 
point of $V$. The point $P$ is a special ramification 
point of type $r+1$ if and only if the largest 
order is at least 
$r+1$, and of type $r-1$ if and only if the largest 
two orders are at least $r-1$.

We say that $P\in C$ is ordinary (resp. a 
{\em Weierstrass point}) if $P$ is an ordinary point 
(resp. a ramification point) of the canonical system, 
i.e. the complete system of sections of $\w_C$.
\end{claim}

\begin{claim}\label{prel2}
({\em Ramification schemes}) 
Let $p\:\TS\to S$ be a smooth, projective 
map of schemes whose fibers are curves. 
For each integer $i\geq 0$, and each 
invertible sheaf $\cal L$ on $\TS$, 
let $J^i_p(\cal L)$ denote the relative 
sheaf of jets, or principal parts, of order $i$ 
of $\cal L$. The sheaf $J^i_p(\cal L)$ is locally free 
of rank $i+1$. Also, there is a natural 
evaluation map, $e_i\:p^*p_*\cal L\to J^i_p(\cal L)$, 
which locally, after trivializations are taken, 
is represented by a Wronskian matrix of functions and 
their derivatives up to order $i$. 

There is a natural 
identification
$J^0_p(\cal L)=\cal L$. 
Furthermore, for each integer $i>0$ there is a 
natural exact sequence of the form:
\begin{equation}\label{trunk}
\begin{CD}
0 \to \w_p^{\otimes i}\otimes\cal L @>>>
J^i_p(\cal L) @>r_i>> J^{i-1}_p(\cal L)\to 0,
\end{CD}
\end{equation}
where $\w_p$ is the relative dualizing sheaf of 
$p$. The truncation maps $r_i$ are compatible with 
the evaluation maps, that is, 
$e_{i-1}=r_i\circ e_i$ for each $i>0$.

Sheaves satisfying the same properties as the 
$J^i_p(\cal L)$ above can be constructed if $p$ is only 
a flat, projective map whose fibers are 
(nodal) curves.
In addition, they 
coincide on the smooth locus of $p$ with the 
corresponding sheaves of jets. These 
sheaves 
appeared in 
\cite{gattothesis}, \cite{estevesthesis}, \cite{lakthor1} 
and \cite{lakthor2}. We will use the same 
notation, $J^i_p(\cal L)$, 
for 
these sheaves.

So, more generally, 
let $p\:\TS\to S$ be a flat, projective map 
whose fibers are curves of genus $g$. 
Let $\cal L$ be an invertible sheaf on $\TS$, and 
$\nu\:\cal V\to p_*\cal L$ any map from a locally free 
sheaf $\cal V$ of constant positive rank, say $r+1$ for 
a certain integer 
$r\geq 0$. For each integer $i\geq 0$, 
consider the natural evaluation map,
$$
\begin{CD}
u_i:p^*\cal V @>p^*\nu >> 
p^*p_*\cal L @>>> J^i_p(\cal L).
\end{CD}
$$
We call the degeneracy scheme of $u_{r+j}$, for $j=-1,1$, 
the {\em special ramification scheme of type 
$r+1+j$} of $(\cal V,\cal L)$, and denote it by 
$VE_j(\cal V,\cal L)$. We call the degeneracy 
scheme of $u_r$ the {\em ramification scheme} of 
$(\cal V,\cal L)$, and denote it by $W(\cal V,\cal L)$.

If $S:=\text{Spec}(\CC)$, if $\TS$ is smooth, and 
if $\nu$ is injective, then the support of 
the scheme 
$W(\cal V,\cal L)$ is the set of ramification points 
of the linear system 
$H^0(S,\cal V)\subseteq H^0(\TS,\cal L)$ of sections 
of $\cal L$. Also, the support of $VE_j(\cal V,\cal L)$ 
is the set of special ramification points of type 
$r+1+j$ of the same linear system, for $j=-1,1$.   

The map $u_r$ 
is a map of locally free sheaves of the same rank $r+1$. 
Taking determinants, $u_r$ induces a 
{\em Wronskian section} $w_p$ of the 
invertible sheaf
$$
\cal W:=\bigwedge^{r+1} J_p^r(\cal L)\otimes
\Big(\bigwedge^{r+1}p^*\cal V\Big)^\vee.
$$
Using the truncation sequences \eqref{trunk}, we get
$$
\cal W\cong\w_p^{\otimes\binom{r+1}{2}}\otimes
\cal L^{\otimes r+1}\otimes
\Big(\bigwedge^{r+1}p^*\cal V\Big)^\vee.
$$
Locally, after trivializations are taken, 
$w_p$ corresponds to a Wronskian determinant of a 
sequence of $r+1$ functions. Its zero scheme is 
the ramification scheme of $(\cal V,\cal L)$. 

The formation of the ramification scheme is functorial 
in the following sense: Suppose there are an invertible 
sheaf $\cal L'$ on $\TS$, a locally free sheaf $\cal V'$ 
of rank $r+1$ on $S$, 
a map $\psi\:\cal L'\to\cal L$, and a commutative 
diagram of maps of the form:
$$
\begin{CD}
\cal V' @>\nu'>> p_*\cal L'\\
@V\mu VV @Vp_*\psi VV\\
\cal V @>\nu >> p_*\cal L.
\end{CD}
$$
Let $V\subseteq S$ 
be the ramification scheme of $\mu$, and 
$Y\subseteq\TS$ that of $\psi$. (So 
$\text{Im}(\psi)=\cal I_{Y|\TS}\otimes\cal L$.) Then
$$
\cal I_{W'|\TS}\cal I_{Y|\TS}^{r+1}=\cal I_{p^{-1}(V)|\TS}
\cal I_{W|\TS},
$$
where $W:=W(\cal V,\cal L)$ and $W':=W(\cal V',\cal L')$.

By derivation, the section $w_p$ induces a global section 
$w'_p$ of the rank-2 locally free sheaf 
$J^1_p(\cal W)$. We will call the zero scheme of this 
section the {\em special ramification scheme} of 
$(\cal V,\cal L)$, and denote it by 
$VSW(\cal V,\cal L)$. 
Notice that the irreducible components 
of the ramification scheme have codimension at most 1 in 
$\TS$, while those of the special ramification schemes 
have codimension at most 2. Also, a local analysis 
of the matrices representing the maps $u_i$ shows that, 
set-theoretically,
$$
VSW(\cal V,\cal L)=VE_{-1}(\cal V,\cal L)\,\bigcup\,
VE_{1}(\cal V,\cal L).
$$

Let $T$ be an $S$-scheme, and let $p_T\:\TS_T\to T$ 
denote the induced famiy by base extension. For each 
coherent sheaf $\cal F$ on $S$ (resp. $\TS$), let 
$\cal F_T$ denote its pullback to $T$ (resp. 
$\TS_T$.) Then $\nu$ induces a map
$$
\nu_T\:\cal V_T\to (p_*\cal L)_T\to p_{T*}\cal L_T,
$$ 
and the (special) ramification scheme(s) of 
$(\cal V,\cal L)$ pull back to the (special) 
ramification scheme(s) of $(\cal V_T,\cal L_T)$. 
Furthermore, if $\cal Y\subseteq\TS_T$ is a $T$-flat 
closed 
subscheme whose fibers over $T$ are subcurves 
of the fibers of $p$, then the (special) ramification 
scheme(s) of $(\cal V,\cal L)$ coincide on 
$\cal Y-(\cal Y\cap\ol{\TS_T-\cal Y})$ with the 
corresponding (special) ramification scheme(s) of 
$(\cal V_T,\cal L_T|_Y)$.

In case $\cal L$ is the relative dualizing sheaf 
of $p$, and $\cal V=p_*\cal L$, the ramification schemes 
and special ramification schemes are called Weierstrass 
schemes and special Weierstrass schemes. In addition, we 
set 
$$
W(p):=W(\cal V,\cal L),\quad VSW(p):=VSW(\cal V,\cal L),
$$
and $VE_j(p):=VE_j(\cal V,\cal L)$ for $j=-1,1$.
\end{claim}

\begin{claim}\label{prel3} ({\em Smoothings}) 
Let $C$ be a curve. A 
{\em smoothing} of $C$ consists of two data: 
a flat, projective map $p\:\TS\to S$ to 
$S:=\text{Spec}(\mathbb C[[t]])$ 
with smooth generic fiber, 
and an isomorphism between the special fiber and $C$. 
The smoothing is called {\em regular} if the total 
space $\TS$ is a regular scheme.

Let $p\:\TS\to S$ be a smoothing of $C$, and identify 
the special fiber with $C$ with the given isomorphism. 
Since the general 
fiber is smooth, for each node $P$ 
of $C$, there are a nonnegative integer $k$ and a 
$\mathbb C[[t]]$-algebra isomorphism
\begin{equation}\label{OP}
\wh{\cal O}_{\TS,P}\cong
\frac{\mathbb C[[t,x,y]]}{(xy-t^{k+1})}.
\end{equation}
We call $k$ the {\em singularity type} 
of $P$ in $\TS$, and 
set $k(P):=k$. Notice that $k(P)=0$ if and only if $\TS$ 
is regular at $P$.

If $E\subseteq C$ is an irreducible component, then 
$E$ is not necessarily a Cartier divisor of $\TS$. 
However, let $m_E$ be the least common 
multiple of the $k(P)+1$ 
for all 
$P\in E\cap\ol{C-E}$. Then there is a natural 
effective Cartier divisor on $\TS$ whose associated 
1-cycle is $m_E[E]$; call 
this divisor $E^p$. At a node 
$P\in E\cap\ol{C-E}$, with $\wh{\cal O}_{\TS,P}$ 
of the form \eqref{OP}, it is given by $x^n=0$, for 
$n:=m_E/(k+1)$, if $x=0$ is 
the equation defining the subcurve $E\subseteq C$. 
Notice that $E$ itself is a Cartier divisor of $\TS$ 
if and only if $m_E=1$.

For each integer $d>0$, let $S\to S$ be the map 
defined by taking $t$ to $t^d$, 
and let $p_d\:\TS_d\to S$ be the smoothing 
induced by base change. 
The special fiber of $p_d$ is equal to that of $p$. 
But, for a node $P\in C$, if $k$ is the 
singularity type of $P$ in $\TS$, then $(k+1)d-1$ is 
the singularity type of $P$ in $\TS_d$.

Suppose $C$ is semistable. 
There are a smoothing $\ol p\:\ol{\TS}\to S$ and an 
$S$-map $b\:\TS\to\ol{\TS}$ that blows down (collapses) 
all rational 
smooth components $E$ of $C$ which intersect $\ol{C-E}$ 
in two points. In fact, just let
$$
\ol{\TS}:=\text{Proj}\Big(\bigoplus_{i\geq 0}
H^0(\TS,\w_p^{\otimes i})\Big),
$$
where $\w_p$ is the relative dualizing sheaf of $p$.
However the singularity types of $\ol{\TS}$ 
are bigger than those of $\TS$: if $\ol C:=b(C)$, and 
$P\in\ol C$ is a node 
obtained by blowing down a chain of $r$ smooth rational 
curves, and $k_0,k_1,\dots,k_r$ are the singularity 
types in $\TS$ of the nodes of $C$ on that chain, 
then the singularity type of $P$ in $\ol{\TS}$ is 
$k_0+k_1+\cdots+k_r+r$.

In certain circumstances, it might be interesting to 
avoid blowing down some of the rational components of 
$C$ in a construction as above. This is possible after 
base change. With a base change we may produce sections 
$\Sigma_i\subset\TS$ of $p$ through its smooth locus 
intersecting the components we do not want to 
blow down. Then just do the above construction with 
$\w_p$ replaced by $\w_p(\sum\Sigma_i)$. 

So, given a node $P$ of $C$, and 
positive integers $m_0,\dots,m_n$, it is possible, 
with base changes, blowups and blowdowns, to 
find an integer $d>0$, a smoothing $\wt p\:\wt\TS\to S$ 
and an $S$-map $b\:\wt\TS\to\TS_d$ such that 
$\wt p=b\circ p_d$ and:
\begin{enumerate}
\item $b$ is an isomorphism off $P$. 
\item $b^{-1}(P)$ is a chain of $n$ smooth rational 
components of the special fiber $\wt C$ of $\wt p$.
\item The singularity types in $\wt\TS$ 
of the nodes $P_0,\dots,P_n$ 
of $\wt C$ on $b^{-1}(P)$, ordered in sequential order, 
are $\ell m_0-1,\,\dots,\,\ell m_n-1$ for a certain 
integer $\ell>0$. (In fact, 
$$
\ell(m_0+\cdots+m_n)=(k+1)d,
$$
where $k$ is the singularity type of $P$ in $\TS$.)  
\end{enumerate} 
\end{claim}

\begin{claim}\label{prel4} ({\em Limit linear series}) 
Let $C$ be a curve, and $p\:\TS\to S$ a smoothing of 
$C$. Identify $C$ with the closed fiber of $p$, and 
denote by $\TS_*$ the general fiber.

Let $\cal L$ be an invertible sheaf on $\TS$. Since 
$p$ is flat, $H^0(\TS,\cal L)$ is a torsion-free 
$\CC[[t]]$-module, whence free. Let 
$V\subseteq H^0(\TS,\cal L)$ be a $\CC[[t]]$-submodule. 
Then also $V$ is free, say of rank $r+1$, for a 
certain integer $r\geq 0$. Assume $V$ is 
\emph{saturated}, i.e. $(V:(t))=V$. Letting 
$V_*$ be the subspace of $H^0(\TS_*,\cal L|_{\TS_*})$ 
generated by $V$, 
we have that $V$ is saturated if and only if 
$V=V_*\cap H^0(\TS,\cal L)$. In our applications 
we will actually have $V=H^0(\TS,\cal L)$, so saturated.

Let $R\subset\TS$ be the ramification scheme of 
$(V\otimes\cal O_S,\cal L)$, as defined in 
Subsection \ref{prel2}. Since 
$\TS_*$ is smooth, $R$ is indeed 
a divisor. However, $R$ may not intersect $C$ 
properly, as $R$ may contain in its support a 
component of $C$. Nevertheless, let 
$\ol R:=\ol{R\cap \TS_*}$. Then $\ol R$ intersects 
$C$ properly. The intersection, 
$\partial R:=\ol R\cap C$ is called the 
\emph{limit ramification scheme}.  

In \cite{EST2} it is shown how to compute the 0-cycle 
$[\partial R]$ associated to $\partial R$ when 
$p$ is regular. 
We review this below.

Let $C_1,\dots,C_n$ be the irreducible components of $C$. 
Since $C$ is connected, for 
each $i=1,\dots,n$ there is an invertible sheaf 
$\cal L_i$ on $\TS$ of the form 
$$
\cal L_i=\cal L\otimes
\cal O_{\TS}(\textstyle\sum_ma_{i,m}C^p_m),
\quad a_{i,m}\in\mathbb Z,
$$
such that the restriction map 
\begin{equation}\label{restri}
H^0(\TS,\cal L_i)\lra H^0(C_i,\cal L_i|_{C_i})
\end{equation}
has kernel $tH^0(\TS,\cal L_i)$. (The divisors $C_m^p$ 
are as explained in Subsection \ref{prel3}.)

There is a natural 
identification $\cal L_i|_{\TS_*}=\cal L|_{\TS_*}$. Using 
it, set
$$
V_i:=H^0(\TS,\cal L_i)\cap V_*\subseteq 
H^0(\TS_*,\cal L|_{\TS_*}).
$$
Then also $V_i$ is saturated and free of rank $r+1$. 
(In fact, $V_{i*}=V_*$.) Let
$\ol V_i\subseteq H^0(C_i,\cal L_i|_{C_i})$ be 
the image of $V_i$ under \eqref{restri}. Since $V_i$ 
is saturated, and \eqref{restri} has kernel 
$tH^0(\TS,\cal L_i)$, the dimension of $\ol V_i$ is 
$r+1$. We call $(\ol V_i,\cal L_i|_{C_i})$ a 
\emph{limit linear system} on $C_i$. 

Let 
$R_i\subseteq C_i$ be the ramification scheme of 
$(\ol V_i,\cal L_i|_{C_i})$, as 
defined in Subsection \ref{prel2}. 
Put $R'_i:=R_i-R_i\cap\ol{C-C_i}$. Then 
$$
[\partial R]\geq [R'_1]+\cdots+[R'_n].
$$
Furthermore, if $p$ is regular, then
\begin{equation}\label{dR}
[\partial R]=\sum_{i=1}^n[R_i]+
\sum_{i<j}\sum_{P\in C_i\cap C_j}(r+1)(r-\ell_{i,j})[P],
\end{equation}
where $\ell_{i,j}:=a_{i,j}+a_{j,i}-a_{i,i}-a_{j,j}$ for 
each distinct $i,j=1,\dots,n$.

When $\cal L=\w$ and $V=H^0(\TS,\w)$, the 
limit ramification scheme is called the 
\emph{limit Weierstrass scheme}, and denoted 
$\partial W_p$; also, a limit linear system is called a 
\emph{limit canonical system}.

Let $P$ be a nonsingular point of $C$, and 
$\Gamma\subset\TS$ a section of $p$ intersecting $C$ at 
$P$. Say, $P\in C_i$. 
Let $P_*$ be the rational point of $\TS_*$ cut out 
by $\Gamma$. Then the behaviour of 
$(V_*,\cal L|_{\TS_*})$ 
at $P_*$ is partially captured by that of 
$(\ol V_i,\cal L_i|_{C_i})$ at $P$. For instance,
we have semicontinuity:
$$
\dim_{\CC}\ol V_i(-aP)\geq\dim_{\CC((t))} V_*(-aP_*)\quad
\text{for each $a=0,1,\dots$.}
$$

In fact, let $m:=\dim_{\CC((t))} V_*(-aP_*)$. Since 
$V_*=V_{i*}$, we may choose a $\CC[[t]]$-basis 
$\sigma_1,\dots,\sigma_m$ of $V_i\cap V_*(-aP_*)$. 
The images $\ol\sigma_i$ in $\ol V_i$ vanish at $P$ with 
multiplicity at least $a$ as well. If there is a nonzero 
$m$-tuple $(c_1,\dots,c_m)\in\CC^m$ such that 
$c_1\ol\sigma_1+\cdots+c_m\ol\sigma_m=0$, then 
\begin{equation}\label{ci}
c_1\sigma_1+\cdots+c_m\sigma_m=t\sigma
\end{equation}
for some 
$\sigma\in V_i$, because \eqref{restri} has kernel 
$tH^0(\TS,\cal L_i)$ and $V_i$ is saturated. Because of 
\eqref{ci}, also $\sigma\in V_i\cap V_*(-aP_*)$, and 
hence $\sigma=b_1\sigma_1+\cdots+b_m\sigma_m$ for certain 
$b_i\in\CC[[t]]$. Plugging this expression in 
\eqref{ci}, we get a nontrivial relation for the sections 
$\sigma_i$, an absurd.
 
In particular, if $P_*$ is a special ramification point 
of type $r+j$ of $V_*$, for $j=-1$ or $j=1$, then so 
is $P$ with respect to $\ol V_i$. When $\cal L=\w$ and 
$V=H^0(\TS,\w)$ we say that $P$ is the limit of a 
special Weierstrass point of type $g+j$ along $p$. 
\end{claim}

\section{General reducible curves}

\begin{proposition}\label{auxprop} 
Let $X$ and $Y$ be two smooth 
nonrational curves. Let $A\in X$ and 
$B\in Y$, and let $C$ be the uninodal curve union 
of $X$ and $Y$ with $A$ identified with $B$. 
Let $p\:\TS\to S$ be a 
smoothing of $C$. Then the following statements hold:
\begin{enumerate}[(i)]
\item\label{auxpropi} 
If $B$ is 
at most a 
simple Weierstrass point
of $Y$, then there is a 
vector
subspace 
$V\subseteq H^0(X,\w_X((g_Y+2)A))$ of 
codimension $1$ containing $H^0(X,\w_X(g_YA))$ such that 
$(V,\w_X((g_Y+2)A))$ is a limit canonical system on $X$. 
Furthermore, 
if $B$ is an ordinary point of $Y$, then 
$A$ is a base point of this system, i.e. 
$V=H^0(X,\w_X((g_Y+1)A))$.
\item\label{auxpropii} 
If $A$ and $B$ are at most simple Weierstrass 
points of $X$ and $Y$, with at least one of them 
ordinary, then the limit Weierstrass scheme 
contains the node of $C$ with multiplicity at most 1.
\end{enumerate}
\end{proposition} 

\begin{proof} Assume $B$ is at most a simple Weierstrass 
point of $Y$. Let $\ell$ be an integer, and set
$\cal L:=\w_p(\ell Y^p)$, 
where $\w_p$ is the relative dualizing sheaf of $p$. 
Then
\begin{equation}\label{LXY}
\cal L|_X\cong\w_X((\ell+1)A)
\quad\text{and}\quad
\cal L|_Y\cong\w_Y((1-\ell)B).
\end{equation}
In addition, the following natural sequences are exact:
	\begin{equation}\label{seq1}
	0\to\cal L|_X(-A)\to\cal L|_C\to\cal L|_Y\to 0,
	\end{equation}
	\begin{equation}\label{seq2}
	0\to\cal L|_Y(-B)\to\cal L|_C\to\cal L|_X\to 0.
	\end{equation}
By Riemann--Roch, $h^0(\cal L|_X)\geq g$ if and only 
if $\ell\geq g_Y$.
On the other hand, by the hypothesis on $B$, we have
$h^0(Y,\cal L|_Y)=\max(0,g_Y+1-\ell)$ for 
$\ell\neq g_Y+1$, whereas for $\ell=g_Y+1$ either 
$h^0(Y,\cal L|_Y)=0$ 
if $B$ is an ordinary point of $Y$, or else 
$h^0(Y,\cal L|_Y)=1$.

Set $\cal M:=\w_p(g_YY^p)$ and 
$\cal N:=\w_p((g_Y+1)Y^p)$. From \eqref{LXY} for 
$\cal L:=\cal M$, and Riemann--Roch, since $g_Y>0$ we 
have $H^1(X,\cal M|_X(-A))=0$. So, the exactness of 
\eqref{seq1} for $\cal L:=\cal M$ implies
        $$h^0(C,\cal M|_C)=h^0(X,\cal M|_X(-A))+
          h^0(Y,\cal M|_Y)=(g-1)+1=g.$$
Thus the restriction 
$H^0(\TS,\cal M)\to H^0(C,\cal M|_C)$ is surjective. 

Consider now the restriction map
        $$\alpha\:H^0(C,\cal M|_C)\longrightarrow
          H^0(X,\w_X((g_Y+1)A)).$$ 
It follows from the exactness of 
\eqref{seq1} that $\alpha$ contains $H^0(X,\w_X(g_YA))$ 
in its image, and from 
the exactness of \eqref{seq2} that the kernel of 
$\alpha$ is isomorphic to $H^0(Y,\w_Y(-g_YB))$. 
Thus $\alpha$ is injective, hence bijective, if and 
only if $B$ is 
an ordinary point of $Y$. In this case, the complete 
linear system of sections of $\w_X((g_Y+1)A)$ is a limit 
canonical system on $X$.
On the other hand, if $B$ is a (simple) 
Weierstrass point of $Y$, the image of 
$\alpha$ is the subspace $H^0(X,\w_X(g_YA))$. 

Applying \eqref{LXY} for $\cal L:=\cal N$, 
as $B$ is at most a simple Weierstrass point of $Y$, 
we get $H^0(Y,\cal N|_Y(-B))=0$. So, the natural map 
        $$\beta\:H^0(C,\cal N|_C)\to 
          H^0(X,\w_X((g_Y+2)A))$$ 
is injective. The maps $\alpha$ and 
$\beta$ fit in a commutative diagram of the form
        $$\begin{CD}
        H^0(\TS,\cal M) @>>> H^0(C,\cal M|_C) @>\alpha >> 
	H^0(X,\w_X((g_Y+1)A))\\
	@VVV @VVV @VVV\\
	H^0(\TS,\cal N) @>>> H^0(C,\cal N|_C) @>\beta >> 
	H^0(X,\w_X((g_Y+2)A)),
	\end{CD}$$
where the horizontal maps are induced by restriction, 
and the vertical maps are induced from the inclusion 
$\cal M\to \cal N$. Since $\beta$ is injective, the 
image $V$ of the bottom composition has codimension 1 
in $H^0(X,\w_X((g_Y+2)A))$, and 
$(V,\w_X((g_Y+2)A))$ is a limit canonical system on $X$. From 
the diagram, $V$ contains the 
image of the top composition, which is 
$H^0(X,\w_X((g_Y+1)A))$ if $B$ is an ordinary point of 
$Y$, and is $H^0(X,\w_X(g_YA))$ otherwise. In the former 
case, by dimension considerations, 
$V=H^0(X,\w_X((g_Y+1)A))$. 
This finishes the proof of \eqref{auxpropi}.

Let us prove \eqref{auxpropii}. Without loss of 
generality, we may assume that $A$ is an 
ordinary point of 
$X$. Then, from Pl\"ucker formula, the ramification 
divisor of the complete linear system of sections 
of $\w_X((g_Y+1)A)$ has degree 
\begin{equation}\label{deg1}
(2g_X+g_Y-1)g+(g_X-1)g(g-1)-g_X
\end{equation}
on $X-A$. On the other hand, since $B$ is at most a 
simple Weierstrass point 
of $Y$, also by Pl\"ucker formula, the ramification 
divisor of the complete linear system of sections 
of $\w_Y((g_X+1)B)$ has degree 
\begin{equation}\label{deg2}
(2g_Y+g_X-1)g+(g_Y-1)g(g-1)-w_B
\end{equation}
on $Y-B$, where $w_B=g_Y$ if $B$ is an ordinary point, 
or else $w_B=g_Y+1$.

Suppose first that $B$ is an ordinary point of 
$Y$. Then, by the already proved \eqref{auxpropi}, 
the limit Weierstrass scheme 
$\partial W_p$ has degree away from the node 
equal to the sum of 
\eqref{deg1} and \eqref{deg2}, 
with
$w_B=g_Y$. 
But this sum is $g^3-g$. So 
the node of $C$ is not contained in $\partial W_p$.

Finally, suppose that $B$ is a (simple) 
Weierstrass point of $Y$. Then, 
by \eqref{auxpropi}, there 
is a vector subspace $V\subset H^0(X,\w_X((g_Y+2)A))$ of 
dimension $g$ containing $H^0(X,\w_X(g_YA))$ such that 
$(V,\w_X((g_Y+2)A))$ is a limit canonical system. 
Since $A$ is an ordinary 
point of $X$, Pl\"ucker formula yields that the 
ramification divisor of $(V,\w_X((g_Y+2)A))$ 
has degree
\begin{equation}\label{deg3}
(2g_X+g_Y)g+(g_X-1)g(g-1)-w_A
\end{equation}
on $X-A$, where $w_A=2g_X+g_Y-1$ if 
$V\neq H^0(X,\w_X((g_Y+1)A))$, 
and $w_A=2g_X+g_Y$ otherwise. At any rate, using 
$w_B=g_Y+1$, the sum of \eqref{deg2} and \eqref{deg3} is 
at least $g^3-g-1$, with equality only if 
$w_A=2g_Y+g_X$. So $\partial W_p$ can only contain the 
node of $C$ with multiplicity 1.
\end{proof}

\begin{lemma}\label{lemaux} 
Let $X$ be a general curve, $P\in X$ a general point, 
and $n$ a positive integer. Let $Q\in X-P$. Then the 
following statements are equivalent, for each $j=-1,1$:
\begin{enumerate}[(i)]
\item\label{lemauxi} The point $Q$ 
is a ramification point 
of the complete system of sections of $\w_X((n+1+j)P)$.
\item\label{lemauxii} There is a unique subspace 
$V\subseteq H^0(X,\w_X((n+2)P))$ of codimension $1$ 
containing $H^0(X,\w_X(nP))$ but different from 
$H^0(X,\w_X((n+1)P))$ such that $Q$ is a special 
ramification point of $(V,\w_X((n+2)P))$ of type 
$g_X+n+j$.
\end{enumerate}
\end{lemma}

\begin{proof} Set $g:=g_X+n$. Also, set 
\begin{equation}\label{V'}
V':=H^0(\w_X(nP))+H^0(\w_X((n+2)P-gQ))\subseteq 
H^0(\w_X((n+2)P)).
\end{equation}
Since $(X,P)$ 
is general, by Prop. 3.1 of \cite{cumestgat1}, 
all the ramification points but $P$ of the 
complete linear system of sections of 
$\w_X(nP)$ or $\w_X((n+2)P)$ are simple. Then
\begin{equation}\label{h00}
h^0(X,\w_X((n+2)P-gQ))=1\,\,\text{ and }\,\,
h^0(X,\w_X(nP-gQ))=0.
\end{equation}
Thus the sum in \eqref{V'} is direct, and $V'$ has 
dimension $g$. In addition,
\begin{equation}\label{V2}
V'(-iQ)=H^0(X,\w_X(nP-iQ))\oplus H^0(X,\w_X((n+2)P-gQ))
\end{equation}
for each $i=0,1,\dots,g$, and thus
\begin{equation}\label{V3}
\dim V'(-iQ)=h^0(X,\w_X(nP-iQ))+1
\quad\text{for each $i=0,1,\dots,g$.}
\end{equation}

Suppose first that \eqref{lemauxi} holds. Then either 
\begin{equation}\label{new1}
h^0(X,\w_X(nP-(g-1)Q))\geq 1,
\end{equation}
in which case \eqref{V3} 
implies that $\dim V'(-(g-1)Q)\geq 2$, and hence 
$Q$ is a special ramification point of 
type $g-1$ of $(V',\w_X((n+2)P))$; or
\begin{equation}\label{new2}
h^0(X,\w_X((n+2)P-(g+1)Q))\geq 1,
\end{equation}
in which case 
\eqref{V'} implies that $V'(-(g+1)Q)\neq 0$, and 
hence $Q$ is a special ramification point of 
type $g+1$ of $(V',\w_X((n+2)P))$. Notice that 
$V'$ cannot be $H^0(X,\w_X((n+1)P))$ because, 
since $(X,P)$ is 
general, the complete linear system of sections of 
$\w_X((n+1)P)$ has no special ramification points 
but $P$, by Prop. 3.1 of \cite{cumestgat1}.

For the uniqueness, just observe that, if $Q$ is a 
ramification point of $(V,\w_X((n+2)P))$, for a subspace 
$V$ as described in \eqref{lemauxii}, then \eqref{h00} 
implies that $V\supseteq H^0(X,\w_X((n+2)P-gQ))$, and 
hence $V\supseteq V'$. Since both $V$ and $V'$ have 
dimension $g$, we have $V=V'$.

Finally, suppose \eqref{lemauxii} holds. As we saw above, 
necessarily $V=V'$. 
So, $Q$ is a special ramification point of 
$(V,\w_X((n+2)P))$ of type $g+j$ if and only if 
$\dim V'(-(g-1)P)\geq 2$ if $j=-1$ or 
$\dim V'(-(g+1)P)\geq 1$ if $j=1$. 
Using \eqref{V3} with $i=g-1$, 
we see that the former inequality 
occurs if and only if \eqref{new1} holds, 
i.e., if and only if $Q$ is a ramification point of the 
complete linear system of sections of $\w_X(nP)$.
On the other hand, since $H^0(X,\w_X(nP-gQ))=0$, 
the latter inequality occurs if and only if 
\eqref{new2} holds, i.e., if and only 
if $Q$ is a ramification point of the 
complete linear system of sections of $\w_X((n+2)P)$.
\end{proof}

\begin{proposition}\label{pro3} Let $Y$ be a 
general smooth curve, $\Delta\subset Y\times Y$ the 
diagonal, and $p_1$ and $p_2$ the projection maps 
from
$Y\times Y$ onto the indicated factors. Set 
$\cal L:=(p_2^*\w_Y)(-(g_Y-1)\Delta)$ and 
$\cal E:=p_{1*}\cal L$. Then $\cal E$ is invertible, 
and the degeneracy scheme of 
the evaluation map $p_1^*\cal E\to \cal L$ 
intersects $\Delta$ transversally along 
the Weierstrass points of $Y$.
\end{proposition}

\begin{proof} By \cite{cumestgat1}, 
Cor.~3.3, the Weierstrass points of 
the general curve are simple. 
Thus $h^0(Y,\w_Y(-(g_Y-1)P))=1$ 
for each $P\in Y$, and hence $\cal E$ is invertible. Let 
$Z$ denote the degeneracy scheme of the evaluation map 
$p_1^*\cal E\to \cal L$. 
For each $P\in Y$, the intersection
$Z\cap p_1^{-1}(P)$ is the ramification scheme of the 
complete linear system of sections of 
$\w_Y(-(g_Y-1)P)$. Thus $Z\cap p_1^{-1}(P)$ is finite 
and contains $(P,P)$ if and only if $P$ is a 
Weierstrass point of $Y$. We need only show now 
that $Z$ intersects $\Delta$ transversally, what will 
follow from showing 
that the intersection number $Z\cdot\Delta$ is 
$g_Y^3-g_Y$. 

Let $\delta\:Y\to Y\times Y$ be the diagonal map. 
We have $\delta^*\cal O_{Y\times Y}(-\Delta)=\w_Y$. Thus 
$$
Z\cdot\Delta=\deg\delta^*Z=
\deg(c_1(\w_Y^{\otimes g_Y})-c_1(\cal E)).
$$
Now, since $Y$ has at most simple Weierstrass points, 
for each $i=0,\dots,g_Y-2$ the natural exact sequence
$$
0\to p_{1*}p_2^*\w_Y(-(i+1)\Delta)\to 
p_{1*}p_2^*\w_Y(-i\Delta)\to 
\w_Y\otimes\delta^*\cal O_{Y\times Y}(-i\Delta)\to 0
$$
is exact. As $c_1(p_{1*}p_2^*\w_Y)=0$ and 
$\delta^*\cal O_{Y\times Y}(-\Delta)=\w_Y$, we get
$$
c_1(\cal E)=
-\big(c_1(\w_Y)+c_1(\w_Y^{\otimes 2})+\cdots+
c_1(\w_Y^{\otimes g_Y-1}\big).
$$
Thus 
$$
Z\cdot\Delta=\sum_{i=1}^{g_Y}i\deg(c_1(\w_Y))=
\binom{g_Y+1}{2}(2g_Y-2)=g_Y^3-g_Y.
$$
\end{proof} 

\begin{theorem}\label{mainthm} Let $X$ and $Y$ be 
two general smooth nonrational curves. Let $A\in X$ and 
$B\in Y$, and let $C$ be the uninodal curve union 
of $X$ and $Y$ with $A$ identified with $B$. 
Set $g:=g_C$. Suppose that 
either $A$ is a general point of $X$ or $B$ is a general 
point of $Y$. Let $Q\in C$ lying on $X$. Then, 
for each $j=-1,1$, the point 
$Q$ is the limit of 
a special Weierstrass point of type 
$g+j$ along a smoothing of $C$ 
if and only if 
$Q$ is not the node of $C$, and either of the following 
two situations occur:
\begin{enumerate}[(i)]
\item\label{mainthmi} 
$Q$ is a special ramification point of 
type $g+j$ of the complete linear system of sections of 
$\w_X((g_Y+1)A)$ or
\item\label{mainthmii} $Q$ is a ramification point of 
the complete linear system of sections of 
$\w_X((g_Y+1+j)A)$ and $B$ is a Weierstrass point of $Y$.
\end{enumerate}
\end{theorem}

\begin{proof} 
We prove the 
``only if'' part of the statement first. Let 
$p\:\TS\to S$ be a smoothing of $C$, as indicated 
in 
Figure~1 
below, 
such that $Q$ is the 
limit of a Weierstrass point of type $g+j$ along $p$. In 
particular, $Q$ appears with multiplicity at least 2 in 
the limit Weierstrass scheme $\partial W_p$.  
\begin{figure}[ht]
\begin{center}
\includegraphics[height=3cm, angle=0]{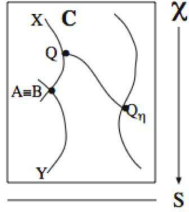}
\caption{The smoothing.}
\end{center}
\end{figure}

Since $X$ and $Y$ are general, their Weierstrass points 
are simple. Also, since either $A$ or $B$ is general, 
either $A$ or $B$ is ordinary. Thus, it follows from 
Proposition \ref{auxprop}, item \eqref{auxpropii}, that 
$Q$ is not the node of $C$.

Suppose first that $B$ is an ordinary point of $Y$. Then, 
by Proposition \ref{auxprop}, item \eqref{auxpropi}, the 
system of sections of $\w_X((g_Y+2)A)$ vanishing at $A$ 
is a limit canonical system, and hence \eqref{mainthmi} 
holds. 

On the other hand, 
suppose that $B$ is a Weierstrass point of $Y$. 
By Proposition \ref{auxprop}, item \eqref{auxpropi}, 
there is a vector subspace 
$V\subset H^0(X,\w_X((g_Y+2)A))$ of codimension 1 
containing $H^0(X,\w_X(g_YA))$ such that 
$(V,\w_X((g_Y+2)A))$ is a limit canonical 
system, and hence admits $Q$ as a special ramification 
point of type $g+j$. Now, 
\eqref{mainthmii} follows from Lemma \ref{lemaux}.

For the ``if'' part of the proof, we will construct 
smoothings as convenient slices of certain $2$-parameter 
families.

Suppose $Q$ is not the node of $C$. Suppose 
first that \eqref{mainthmi} holds. 
Then $g_X\geq 2$. Also, 
it follows from 
Prop. 3.1 in~\cite{cumestgat1} that $A$ is not a general 
point of $X$. So, by hypothesis, $B$ is a general point 
of $Y$, whence an ordinary point. 

We will first deform $C$ by letting $A$ vary to a general 
point. More precisely, let $\Delta\subseteq X\times X$ 
be the diagonal, 
and consider the union $U$ of $X\times X$ with 
$Y\times X$ with $\Delta$ naturally identified 
with $\{B\}\times X$. Let $q\:U\to X$ be the 
projection onto the second factor. 
Since $X$ is nonsingular, we may identify 
the complete local ring of $X$ at $A$ 
with $\CC[[t]]$, and let 
$\wt q\:\wt U\to S$ be the family induced over 
$S:=\text{Spec}(\CC[[t]])$ by base change.

Let $V:=\CC[[t_1,t_2,\dots,t_N]]$ be 
the base of the universal 
deformation space of $C$, where $t_1=0$ 
corresponds to equisingular 
deformations. The map $\wt q$ 
corresponds to a local homomorphism 
$h\:V\to\CC[[t]]$ such that $h(t_1)=0$. 
Since $g_X\geq 2$, the map $\wt q$ is not, even 
infinitesimally, a constant family. So 
there is $j\geq 2$ such that 
$h(t_j)$ generates $t\CC[[t]]$. 
We may assume that $j=2$ and, after a 
harmless reparameterization, that $h(t_2)=t$. 
Letting $p_i(t):=h(t_i)$ 
for each $i\geq 3$, we have 
$h(t_i-p_i(t_2))=0$ for each $i\geq 3$. 
Consider the two-parameter 
subfamily of the universal deformation of $C$ 
given precisely by the equations 
$t_i-p_i(t_2)=0$ for $i=3,\dots,N$. 
Identify the base of this family 
with $S_2:=\text{Spec}(\CC[[t_1,t_2]])$, 
and let $u\:T\to S_2$ denote the 
map giving the family, which is depicted 
in 
Figure~2 
below. 
\begin{figure}[ht]
\begin{center}
\includegraphics[height=4.5cm, angle=0]{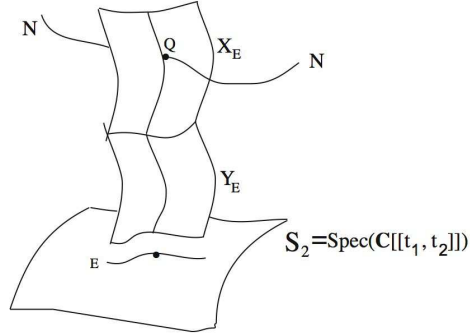}
\caption{The first family.}
\end{center}
\end{figure}

Notice that $T$ is a regular threefold. Let 
$E\subset S_2$ be the 
Cartier divisor given by $t_1=0$. The slice 
$u_E\: u^{-1}(E)\to E$ is precisely $\wt q$, 
under the identification $t_2=t$. 
Hence, the pullback 
$\pi^*E$ is the sum of two 
effective Cartier divisors, 
$X_E$ and $Y_E$, the first isomorphic to 
$X\times E$, the second to 
$Y\times E$, whose intersection on 
$Y\times E$ is $B\times E$, and 
on $X\times E$ is the graph $\Sigma$ 
of a nonconstant map 
$E\to X$, sending the special point 
$o\in E$ to $A$, and the 
general point $e\in E$ to 
the general point of $X$. 

Let $\cal M:=\w_u(g_YY_E)$, where $\w_u$ 
is the relative dualizing sheaf 
of $u$. 
Then
\begin{equation}\label{MXE}
\cal M|_{X_E}\cong\w_{X_E/E}((g_Y+1)\Sigma).
\end{equation}
A fiberwise check, as done in the proof of 
Proposition~\ref{auxprop}, shows 
that $u_*\cal M$ is locally free of 
rank $g$, with formation 
commuting with base change. In addition, since $B$ is 
an ordinary point of $Y$, the natural map 
$$
\gamma\:(u_*\cal M)|_E\longrightarrow 
u_{E*}(\cal M|_{X_E})
$$ 
is an isomorphism. 

Form the special ramification scheme 
$Z\subseteq T$ 
of type $g+j$ of 
$(u_*\cal M,\cal M)$, 
as explained in 
Subsection~\ref{prel2}. Since 
$\gamma$ is an isomorphism, $Z$ agrees on 
$X_E-\Sigma$ with the special 
ramification scheme of 
type $g+j$ of 
$(u_{E*}(\cal M|_{X_E}),\cal M|_{X_E})$. 
Because of \eqref{MXE}, the fact that 
$\Sigma\cap u^{-1}(o)=\{A\}$, and the hypothesis on $Q$, 
we have that 
$Q$ is an isolated point of $Z\cap u^{-1}(o)$. 
Furthermore,
since the general point of 
$\Sigma$ is the general point of 
$X\times\{e\}$, Prop. 3.1 in 
\cite{cumestgat1} yields 
$Z\cap u^{-1}(e)\subseteq Y\times\{e\}$. 

Since 
$Z$ is defined locally by two regular functions, there 
is an irreducible curve $N\subseteq Z$ 
containing $Q$. Since 
$Q$ is an isolated point of $Z\cap u^{-1}(o)$, and 
$Z\cap u^{-1}(e)\subseteq Y\times\{e\}$, 
the general point of 
$N$ must lie on a smooth fiber of 
$u$, and hence be a special 
Weierstrass point of type $g+j$ of that fiber. 
So $Q$ is the limit of a special 
Weierstrass point of type $g+j$.

Suppose now that \eqref{mainthmii} holds. In particular, 
$B$ is a Weierstrass point of 
$Y$, and hence $g_Y\geq 2$. Letting $B$ vary, we may
construct a 
two-parameter family similar to the 
one constructed in the first case. 
Thus we get a family of curves 
$u\:T\to S_2$ over 
$S_2=\text{\rm Spec}(\CC[[t_1,t_2]])$ such that 
$T$ is a regular threefold, 
and the pullback $u^*E$ of the Cartier 
divisor $E\subset S_2$ given by $t_1=0$ is the sum 
of two effective Cartier divisors, 
$X_E$ and $Y_E$, the first 
isomorphic to $X\times E$, the 
second to $Y\times E$, whose 
intersection on $X\times E$ is 
$A\times E$, and on $Y\times E$ is 
the graph $\Sigma$ of a nonconstant 
map $E\to Y$, sending the 
special point $o\in E$ to $B$, 
and the general point $e\in E$ to the 
general point of $Y$.

Let $\wt S_2\to S_2$ be the blowup 
map of $S_2$ at $o$, 
and denote by $F\subset\wt S_2$ the 
exceptional divisor. Abusing 
notation, we denote 
the strict transform of $E$ by $E$ as well, 
and let $o$ denote the point of 
intersection of $E$ and $F$. The fibered product 
$T\times_{S_2}\wt S_2$ is singular 
only at the node of the fiber over $o$ 
of the second projection 
$T\times_{S_2}\wt S_2\to\wt S_2$.

Let $\wt T$ be the blowup of 
$T\times_{S_2}\wt S_2$ along the subscheme 
$Y_E\times_{S_2}\wt S_2\subset T\times_{S_2}\wt S_2$. 
A local analysis 
shows that $\wt T$ is smooth. 
Denote by $\wt X_E$ and $\wt Y_E$ the 
strict transforms 
in $\wt T$
of $X_E\times_{S_2}E$ and 
$Y_E\times_{S_2}E$. Let also 
$\wt X_F$ and $\wt Y_F$ denote 
the strict transforms of $X\times F$ 
and $Y\times F$. Let $\wt u\:\wt T\to\wt S_2$ be the 
induced map. The fiber $\wt T_o:=\wt u^{-1}(o)$ consists 
of three components: two of 
them disjoint and naturally identified 
with $X$ and $Y$, and the 
remaining, say $R$, isomorphic 
to a line and meeting 
$X$ and $Y$ transversally at $A$ 
and $B$. A local analysis shows that 
$\wt X_E\cap\wt T_o=X$ and $\wt Y_E\cap\wt T_o=Y\cup R$,
while 
$\wt X_F\cap\wt T_o=X\cup R$ and $\wt Y_F\cap\wt T_o=Y$.
Figure~3
below describes the family 
given by $\wt u$.
\begin{figure}[ht]
\begin{center}
\includegraphics[height=6cm, angle=0]{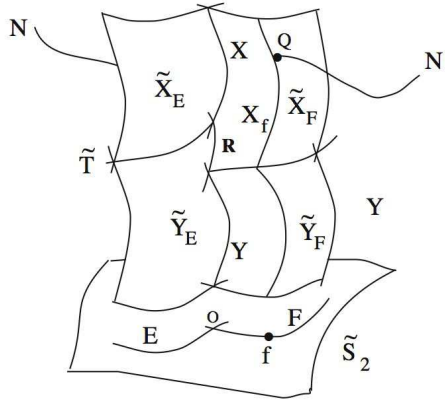}
\caption{The second family.}
\end{center}
\end{figure}
   
For each $z\in E\cup F$, 
let $X_z$ and $Y_z$ denote the 
components of $\wt u^{-1}(z)$ 
that are base extensions of $X$ and $Y$.

Let $\w_{\wt u}$ 
be the relative dualizing sheaf of $\wt u$. Let 
$$
\cal M:=\w_{\wt u}(g_Y(\wt Y_E+\wt Y_F)),\quad
\cal N:=\cal M(\wt Y_F),\quad
\text{and}
\quad\cal P:=\cal N(\wt Y_E).
$$ 
Clearly, $\cal M\subset\cal N\subset\cal P$. 
The restriction $\cal P|_{\wt X_F}$ is the 
pullback of $\w_X((g_Y+2)A)$ under the composition 
$\wt X_F\to X\times F\to X$. Thus 
$$
\wt u_*(\cal P|_{\wt X_F})=
H^0(X,\w_X((g_Y+2)A)\otimes\cal O_F,
$$
and in particular 
$\wt u_*(\cal P|_{\wt X_F})$ is a locally free 
$\cal O_F$-module of rank $g+1$.

We claim that the natural composition
$$
\delta\:(\wt u_*\cal N)|_F\longrightarrow
(\wt u_*\cal P)|_F\longrightarrow
\wt u_*(\cal P|_{\wt X_F})
$$ 
is injective with invertible cokernel, 
and that, as $f$ ranges in $F-o$, 
the image $V_f$ of $\delta(f)$ ranges through 
all subspaces of $H^0(X,\w_X((g_Y+2)A))$ of dimension $g$ 
containing $H^0(X,\w_X(g_YA))$ but 
$H^0(X,\w_X((g_Y+1)A))$. In particular, 
$(\wt u_*\cal N)|_F$ is locally free of rank $g$.

Once the claim is established, 
we 
proceed as
in the first case. Indeed, a 
fiberwise analysis shows that 
$\wt u_*\cal N$ is locally free 
of rank $g$ on $\wt S_2-F$. 
Thus, from the claim, $\wt u_*\cal N$ 
is locally free of rank $g$ 
everywhere. Form the special ramification scheme $Z$ of 
type $g+j$ of $(\wt u_*\cal N,\cal N)$. 
For each 
$f\in F-o$, since $\delta(f)$ is injective, $Z$ agrees 
on $X_f-A$ with the special ramification scheme of 
type $g+j$ of $(V_f,\w_X((g_Y+2)A))$.

Now, by Lemma \ref{lemaux}, there is a subspace 
$V\subseteq H^0(X,\w_X((g_Y+2)A))$ of codimension $1$ 
with $V\supset H^0(X,\w_X(g_YA))$ but 
$V\neq H^0(X,\w_X((g_Y+1)A))$ such that $Q$ is a special 
ramification point of $(V,\w_X((g_Y+2)A))$ of type 
$g+j$. From the claim there is $f\in F-o$ 
such that $V=V_f$. So, viewing $Q$ as a point 
of 
$X_f$,
we have $Q\in Z$. 

Since all irreducible components of $Z$ have 
codimension at most 2 in $\wt T$, there is an 
irreducible curve $N\subseteq Z$ 
passing through $Q\in X_f$. 
Now, only finitely many points of $X$
can be special ramification points of type $g+j$ 
of a linear system like $V$, 
namely, by Lemma \ref{lemaux}, 
the ramification points of the complete linear system of 
sections $\w_X((g_Y+1+j)A)$. But, 
again by Lemma~\ref{lemaux}, 
each of these points is a special ramification 
point of a unique $V$. Thus, 
for all but finitely many $f\in F-o$, 
the image $V_f$ has no 
special ramification points of type $g+j$. Hence 
$N$ intersects only finitely many 
fibers of $\wt u$ over $F$. So 
the general point of $N$ must 
be on a smooth fiber of $\wt u$, 
and hence be a special Weierstrass 
point of type $g+j$ of that fiber. So $Q$ is 
the limit of a special Weierstrass point of type $g+j$.

Now, let us establish the claim. First, 
a fiberwise analysis shows that 
$\wt u_*\cal M$ is locally free of rank $g$, and that 
$R^1\wt u_*\cal M$ is invertible, both with formation 
commuting with base change. 
Consider the long exact sequence 
in higher direct images:
	$$0\to\wt u_*\cal M\to\wt u_*\cal N\to
        \wt u_*(\cal N|_{\wt Y_F})\to R^1\wt u_*\cal M
	\to R^1\wt u_*\cal N\to 
	R^1\wt u_*(\cal N|_{\wt Y_F})\to 0.$$
Since $R^1\wt u_*\cal M$ is invertible, and 
$\wt u_*(\cal N|_{\wt Y_F})$ is supported on 
$F$, the middle map 
above is zero, breaking up the long 
sequence in two short exact 
sequences,
	$$0\to\wt u_*\cal M\to\wt u_*\cal N\to
	\wt u_*(\cal N|_{\wt Y_F})\to 0,$$
	$$0\to R^1\wt u_*\cal M\to R^1\wt u_*\cal N\to 
	R^1\wt u_*(\cal N|_{\wt Y_F})\to 0.$$
The exactness of the first sequence shows 
the surjectivity of the natural map 
$(\wt u_*\cal N)|_F\to\wt u_*(\cal N|_{\wt Y_F})$. 
Now, a fiberwise analysis, using that $B$ is a simple 
Weierstrass point of $Y$, shows that 
$R^1\wt u_*(\cal N|_{\wt Y_F})$ is a 
locally free $\cal O_F$-module of rank 2. 
So, since $R^1\wt u_*\cal M$ is also locally free, 
the exactness of the second sequence above implies 
that, for each Cartier divisor $G\subset\wt S_2$ 
intersecting $F$ properly, 
the natural multiplication-by-$G$ map 
$(R^1\wt u_*\cal N)(-G)\to R^1\wt u_*\cal N$ 
is injective, and hence the natural map 
        $$\delta_G\:\wt u_*\cal N|_G\longrightarrow
        \wt u_*(\cal N|_{\wt u^{-1}(G)})$$ 
is an isomorphism. This isomorphism 
allows us to work with slices of 
the family $\wt u$ that intersect $F$ properly. 

In particular, for 
each $f\in F-\{o\}$, let $G\subset\wt S_2$ be a smooth curve 
passing through $f$, 
and whose general point lies on 
$\wt S_2-(E\cup F)$. So we have a 
smoothing $\wt u_G\:\wt u^{-1}(G)\to G$ of the fiber 
$C$, and we can also choose $G$ such that 
$\wt u_G$ is regular. Then, 
as we saw in the proof of Proposition \ref{auxprop}, 
the natural map
        $$\delta_{G,f}\:(\wt u_{G*}
          (\cal N|_{\wt u^{-1}(G)}))(f)\to 
          H^0(X_f,\cal N|_{X_f})$$
is injective and, under the isomorphism 
$\cal N|_{X_f}\cong\w_X((g_Y+2)A)$, 
its image is a $g$-dimensional subspace of 
$H^0(X,\w_X((g_Y+2)A))$ 
that contains $H^0(X,\w_X(g_YA))$ but is different from 
$H^0(X,\w_X((g_Y+1)A))$. Now, since 
$\delta(f)=\delta_{G,f}\circ\delta_G(f)$, and 
$\delta_G$ is an isomorphism, $\delta(f)$ 
has the same properties.

To understand what happens at $o$, 
consider the slice of $\wt u$ 
over $E$. We claim the natural map 
$$
\eta\:\wt u_*(\cal N|_{\wt u^{-1}(E)})\longrightarrow
\wt u_*(\cal N|_{\wt X_E})
$$ 
is an isomorphism. Indeed, let 
$\Sigma_E:=\wt X_E\cap\wt Y_E$. 
Since $\delta_E$ is an isomorphism, applying 
the long exact sequence in 
higher direct images to the exact sequence
	$$0\to\cal N|_{\wt X_E}(-\Sigma_E)\to
          \cal N|_{\wt u^{-1}(E)}
          \to\cal N|_{\wt Y_E}\to 0$$
we get that the natural map 
$(\wt u_*\cal N)|_E\to\wt u_*(\cal N|_{\wt Y_E})$ 
is surjective, and 
that the image of $\eta$ contains 
$\wt u_*(\cal N|_{\wt X_E}(-\Sigma_E))$. 
Thus, to show our last claim 
we need only show that the natural map 
$$
\epsilon\:\wt u_*(\cal N|_{\wt Y_E})\longrightarrow
\wt u_*(\cal N|_{\Sigma_E})
$$
is an isomorphism.

Since the map $\Sigma_E\to E$ is an isomorphism, 
$\wt u_*(\cal N|_{\Sigma_E})$ is locally free of rank 1. 
Also $\wt u_*(\cal N|_{\wt Y_E})$ is locally free of 
rank 1, because it is so over the generic point $e\in E$. 
Since the point in 
the intersection $\Sigma_E\cap\wt u^{-1}(o)$
is not a Weierstrass point 
of
$Y_e$, the map $\epsilon(e)$ is an isomorphism. 

To show that 
also
$\epsilon(o)$ is an 
isomorphism, it amounts 
to show that 
the point in 
$\Sigma_o:=\Sigma_E\cap\wt u^{-1}(o)$ 
of 
intersection of $X_o$ and $R$ is 
not a limit ramification point of 
$(\wt u_*\cal N|_{\wt Y_E},\cal N|_{\wt Y_E})$.
This is indeed the case, since $\wt Y_E$ 
is the blowup of $Y\times E$ 
at $(B,o)$, and $\Sigma_E$ is the 
strict transform of the graph 
of the map $E\to Y$ obtained by 
considering the identity map of $Y$ 
locally analytically at $B$. So, the transversality 
stated in Proposition \ref{pro3} 
shows 
that $\Sigma_o$ is not a limit ramification point.

Finally, 
since $\eta$ and $\delta_E$ are isomorphisms, it follows 
that $\delta(o)$, which is the composition of 
the isomorphism $\eta(o)\circ\delta_E(o)$ with the 
inclusion 
$$
H^0(X,\w_X((g_Y+1)A))\to H^0(X,\w_X(g_Y+2)A),
$$ 
is 
injective of rank $g$. So, $\delta(f)$ is injective 
of rank $g$ for every $f\in F$, and hence $\delta$ is 
injective with invertible cokernel. Moreover, as 
the image of $\delta(o)$ 
is different from that of $\delta(f)$ for $f\in F-o$, 
then, as $f$ varies in $F-o$,
the image $V_f$ of $\delta(f)$ varies 
through all the 
codimension-1 
subspaces of $H^0(X,\w_X((g_Y+2)A))$ 
containing 
$H^0(X,\w_X(g_YA))$, with the exception of 
$H^0(X,\w_X((g_Y+1)A))$. The 
proof of the 
claim
is complete.
\end{proof}

\section{General irreducible singular curves}

\begin{proposition}\label{2pt} Let $a$ and $b$ be 
positive integers. Let $X$ be a general 
smooth curve of genus $g\geq 0$, and $P$ and $Q$ general 
points on $X$. Then the 
complete linear system of sections of $\w_X(aP+bQ)$ 
has only simple ramification points, and $P$ and $Q$ 
are not among them.
\end{proposition}

\begin{proof} If $g=0$, all complete linear systems on 
$X$ have no ramification points. If $g=1$, 
the curve $X$ can be any curve of genus 1, as long as 
$P-Q$ is neither $a$-torsion nor $b$-torsion in 
the Jacobian variety. 

Assume $g>1$. Let $i<g$ be any positive integer, and 
put $j:=g-i$. 
Let $Y$ and $Z$ be two general 
smooth curves, $Y$ of genus 
$i$, and $Z$ of genus $j$, and let $A$ and $M$ be 
general points of $Y$, and $B$ and $N$ general points 
of $Z$. By induction 
on the genus, we may assume 
the statement of the proposition 
holds for $(Y,A,M)$ and $\w_Y(aA+(b+j)M)$, and 
for $(Z,B,N)$ and $\w_Z(bB+(a+i)N)$.

Let $X_0$ be the nodal curve of genus $g$ given as 
the union of $Y$ and $Z$, with $M$ identified with 
$N$. Since $X_0$ is nodal, and $A$ and $B$ are 
nonsingular points of $X_0$, there are a 
regular smoothing $p\:\TS\to S$ of $X_0$, and 
sections $\Gamma,\,\Delta\subset\TS$ such that, 
identifying the closed fiber of $p$ with $X_0$, we 
have $\Gamma\cap X_0=\{A\}$ and $\Delta\cap X_0=\{B\}$. 

Let $\TS_*$ denote the general fiber of $p$. 
Let $P$ and $Q$ be the 
points of intersection of $\Gamma$ and $\Delta$ with 
$\TS_*$. 
The 2-pointed curve $(\TS_*,P,Q)$ is 
defined over a finitely generated field extension 
of $\QQ$, and hence can be viewed as a 2-pointed 
complex curve. We claim the statement of the 
proposition holds for this two-pointed curve.  

To prove our claim, let 
$\w_p$ be the relative dualizing sheaf of 
$p$. Let $W_*\subset\TS_*$ be the ramification 
divisor of the complete linear system 
of sections of $\w_p(a\Gamma+b\Delta)|_{\TS_*}$. 
We need only show that $W_*$ is reduced, and does not 
contain $P$ or $Q$ in its support. For this, it is enough 
to show that the limit Weierstrass scheme 
$\partial W$ is reduced and 
does not contain $A$ or $B$ in its support. 

Since $\TS$ is regular, $Y$ and $Z$ are Cartier 
divisors. Set
\begin{eqnarray*}
{\cal L}_1&:=&\w_p(a\Gamma+b\Delta+(b+j-1)Z),\\
{\cal L}_2&:=&\w_p(a\Gamma+b\Delta+(a+i-1)Y).
\end{eqnarray*}
Then
\begin{eqnarray*}
{\cal L}_1|_Y=\w_Y(aA+(b+j)M),&\quad&
{\cal L}_1|_Z=\w_Z(bB+(2-b-j)N),\\
{\cal L}_2|_Z=\w_Z(bB+(a+i)N),&\quad&
{\cal L}_2|_Y=\w_Y(aA+(2-a-i)M).
\end{eqnarray*}
Due to the generality of $M$ and $N$, we have
$$
h^0(Y,{\cal L_2}|_Y(-M))=h^0(Z,{\cal L}_1|_Z(-N))=0,
$$
and hence the natural maps
$$
\frac{H^0(\TS,{\cal L}_1)}{tH^0(\TS,{\cal L}_1)}\lra 
H^0(Y,{\cal L}_1|_Y)\quad{\rm and}\quad
\frac{H^0(\TS,{\cal L}_2)}{tH^0(\TS,{\cal L}_2)}\lra 
H^0(Y,{\cal L}_2|_Z)
$$
are injective. They are actually isomorphisms, 
since $H^0(\TS,{\cal L}_i)$ is free of rank 
$g+a+b-1$, and
$$
h^0(Y,{\cal L_1}|_Y)=h^0(Z,{\cal L}_2|_Z)=g+a+b-1,
$$
by the Riemann--Roch theorem. 

Since
$$
{\cal L}_2\cong{\cal L}_1((g+a+b-2)Y),
$$
it follows that $[\partial W]=[R_1]+[R_2]$, where 
$R_1$, resp. $R_2$, 
is the ramification divisor of the complete 
linear system of sections of $\w_Y(aA+(b+j)M)$, 
resp. $\w_Z(bB+(a+i)N)$; see Subsection \ref{prel4}. 
By induction, 
$R_1$ and $R_2$ are reduced and, viewed as subschemes 
of $X_0$, disjoint. So $\partial W$ is reduced. 
In addition, $A$ and $B$ are not in the supports of 
$R_1$ and $R_2$, and thus are not in support of 
$\partial W$ either.
\end{proof}

\begin{theorem}\label{irredlim} 
Let $C$ be a nodal curve of genus $g\geq 2$ 
that is union of a smooth curve $X$ of genus $g-1$ and a 
chain $E=(E_1,\dots,E_{g-1})$ of $g-1$ rational curves. 
Suppose $X$ meets $E$ only in $E_1$ and $E_{g-1}$, 
at points 
$A$ and $B$, respectively. Suppose $X$ is general, and 
$A$ and $B$ are general points of $X$. 
Let $p:\TS\to S$ be a regular smoothing 
of $C$ and $W(p)$ its Weierstrass scheme. 
Then $W(p)$ is a Cartier divisor, and the difference
$$
W(p)-\sum_{i=1}^{g-1}\frac{i(g-i)g}{2}E_i
$$ 
is effective and intersects each fiber 
of $p$ transversally. 
In particular, the limit Weierstrass scheme 
of $p$ is reduced 
and contains no node of $C$ in its support.
\end{theorem}

\begin{proof} 
Figure~4 
below describes the 
curve $C$.
\begin{figure}[ht]
\begin{center}
\includegraphics[height=3cm, angle=0]{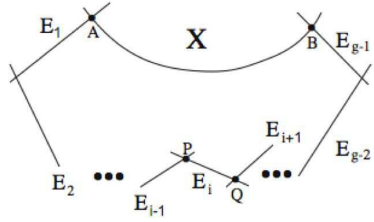}
\caption{The nodal curve.}
\end{center}
\end{figure}

We will prove the second statement 
first. Let $\w_p$ be the relative dualizing sheaf 
of $p$. By adjunction, $\w_p$ restricts to the 
trivial sheaf 
on each $E_i$ and to $\w_X(A+B)$ on $X$. 
So, each global section 
of $\w_p$ that vanishes on $X$ vanishes on 
the whole fiber 
$C$, and hence the restriction map
$$
H^0(\TS,\w_p)\to H^0(X,\w_X(A+B)) 
$$
has image of dimension $g$. Since $h^0(X,\w_X(A+B))=g$ 
by the Riemann--Roch theorem, the restriction map is 
surjective. So the complete linear system of 
sections of $\w_X(A+B)$ is a limit canonical system on 
$X$.

By Proposition \ref{2pt}, the ramification points of 
$H^0(X,\w_X(A+B))$ 
are simple. So, the limit Weierstrass scheme 
$\partial W$ of $f$ is reduced 
on $X-\{A,B\}$. Also, since neither $A$ nor $B$ is 
a ramification point of $H^0(X,\w_X(A+B))$, 
again by Proposition \ref{2pt}, Pl\"ucker formula yields
\begin{equation}\label{WfX}
\deg(\partial W\cap X-\{A,B\})=g^3-g^2.
\end{equation}

For each $i=1,\dots,g-1$ set
\begin{equation}\label{Li}
\cal L_i:=\w_p\Big(-i(g-i)E_i-\sum_{m=1}^{i-1}m(g-i)E_m
-\sum_{m=i+1}^{g-1}(g-m)iE_m\Big).
\end{equation}
Notice that $\cal L_i$ has degree $g$ on $E_i$, 
zero on each $E_m$ for $m\neq i$, and 
$$
\cal L_i|_X\cong\w_X\big(-(g-i-1)A-(i-1)B\big).
$$
Since $A$ and $B$ are general 
points of $X$, we have $h^0(X,\cal L_i|_X)=1$ and 
\begin{equation}\label{-A-B}
h^0\big(X,\cal L_i|_X(-A)\big)=
h^0\big(X,\cal L_i|_X(-B)\big)=
h^0\big(X,\cal L_i|_X(-A-B)\big)=0.
\end{equation}
Let $V_i$ be the image of 
$$
H^0(\TS,\cal L_i) \lra H^0(E_i,\cal L_i|_{E_i})
$$
Since $h^0(X,\cal L_i|_X(-A-B))=0$, and 
$\deg\cal L_i|_{E_m}=0$ for $m\neq i$, 
each global section of 
$\cal L_i$ that vanishes on $E_i$ 
vanishes on the whole $C$. 
Thus $(V_i,\cal L_i|_{E_i})$ is a limit canonical 
system on $E_i$.

Let $P$ and $Q$ be the nodes of $C$ in $E_i$. 
A section in $V_i$ that vanishes at $P$ (or $Q$) 
is the restriction of a global 
section of $\cal L_i$ that vanishes at $P$ (or $Q$). 
Since 
$\deg\cal L_i|_{E_m}=0$ for $m\neq i$, it follows from 
\eqref{-A-B} that this global section vanishes on all 
components of $C$ but possibly $E_i$. In particular, its 
restriction in $V_i$ vanishes at $P$ and $Q$.
In other words,
\begin{equation}\label{ViPi}
V_i(-P)=V_i(-Q)=V_i(-P-Q)=
H^0(E_i,\cal L_i|_{E_i}(-P-Q)),
\end{equation}
where the last equality follows from 
dimension considerations. 

By \eqref{ViPi}, the system $V_i$ contains a complete 
subsystem of codimension 1, 
namely $H^0(E_i,\cal L_i|_{E_i}(-P-Q))$. Since 
complete systems on the projective 
line have no ramification 
points, the order sequence of $V_i$ at each point of 
$E_i$ starts with $0,\dots,g-2$. 
The last order can only be 
$g-1$ or $g$, since $\cal L_i|_{E_i}$ 
has degree $g$. Thus all 
ramification points of $V_i$ are simple. 
In addition, since 
$P$ and $Q$ are not ramification points of 
$H^0(E_i,\cal L_i|_{E_i}(-P-Q))$, it follows from 
\eqref{ViPi} that they are not ramification points of 
$(V_i,\cal L_i|_{E_i})$ either.

Thus $\partial W$ is 
reduced on $E_i-\{P,Q\}$ and, since neither 
$P$ nor $Q$ is a ramification point of 
$(V_i,\cal L_i|_{E_i})$, Pl\"ucker formula yields
\begin{equation}\label{WfEi}
\deg(\partial W\cap E_i-\{P,Q\})=g.
\end{equation}

Finally, since there are 
$g-1$ rational components in $C$, 
we get
$$
g^3-g=\deg(\partial W\cap C)\geq g^3-g^2+(g-1)g=g^3-g,
$$
where the inequality 
follows from combining \eqref{WfX} and 
\eqref{WfEi} for $i=1,\dots,g-1$. 
The inequality is thus an 
equality, showing that we accounted for all points in 
the support of $\partial W$, and hence that 
all of them appear with multiplicity 1. 
The second statement is proved.

To prove the first statement, we will consider a 
filtration $\cal L_{i,j}$ of subsheaves of $\w_p$ 
containing $\cal L_i$, defined below.

For each $i=1,\dots,g-1$ and each 
$j=0,1,\dots,i(g-i)-1$, let $k,\,k',\,\ell,\,\ell'$ be 
integers such that
\begin{equation}\label{jkl}
j=ki+\ell=k'(g-i)+\ell',\quad 
0\leq k,\ell'\leq g-i-1,\quad
0\leq k',\ell\leq i-1,
\end{equation}
and put
\begin{eqnarray}
&&c_{i,j,m}:=km+\max(0,\ell-i+m+1),\quad
m=1,\dots,i,\label{cijm1}\\
&&c'_{i,j,m}:=k'(g-m)
+\max(0,\ell'+i-m+1),\quad m=i,\dots,g-1.
\label{cijm2}
\end{eqnarray}
Notice that $c_{i,j,i}=c'_{i,j,i}=j+1$. Finally, set
$$
D_{i,j}:=\sum_{m=1}^ic_{i,j,m}E_m+
\sum_{m=i+1}^{g-1}c'_{i,j,m}E_m,
$$
and put $\cal L_{i,j}:=\w_f(-D_{i,j})$. Notice that 
$$
\cal L_i=\cal L_{i,i(g-i)-1}.
$$

For each $i=1,\dots,g-1$ set $D_{i,-1}:=0$ and 
$\cal L_{i,-1}:=\w_f$. And for each 
$j=0,1,\dots,i(g-i)-1$ set $F_{i,j}:=D_{i,j}-D_{i,j-1}$. 
Then $\cal L_{i,j}=\cal L_{i,j-1}(-F_{i,j})$.  
It follows from \eqref{jkl}, \eqref{cijm1} and 
\eqref{cijm2} that
$$
F_{i,j}=E_{i-\ell}+E_{i-\ell+1}+\cdots+E_i
+E_{i+1}+\cdots+E_{i+\ell'}.
$$
Using this, it can also be shown, by induction on $j$, 
that 
$$
\cal L_{i,j}|_X\cong\begin{cases}
\w_X((1-k)A+(1-k')B)&\text{if $\ell\neq i-1$ and 
$\ell'\neq g-i-1$},\\
\w_X(-kA+(1-k')B)&\text{if $\ell=i-1$ and 
$\ell'\neq g-i-1$},\\
\w_X((1-k)A-k'B)&\text{if $\ell\neq i-1$ and 
$\ell'=g-i-1$},\\
\w_X(-kA-k'B)&\text{if $\ell=i-1$ and 
$\ell'=g-i-1$}.\\
\end{cases}
$$
and
$$
\deg\cal L_{i,j}|_{E_m}=\begin{cases}
k+k'+2&\text{if $m=i$},\\
-1&\text{if $m=i-\ell-1$ or $m=i+\ell'+1$},\\
0&\text{otherwise}.
\end{cases}
$$

It follows that
\begin{equation}\label{h0F}
h^0(F_{i,j},\cal L_{i,j-1}|_{F_{i,j}})=k+k'+1,
\end{equation}
and, setting $\wh F_{i,j}:=\ol{C-F_{i,j}}$, 
$$
h^0\big(\wh F_{i,j},\cal L_{i,j}|_{\wh F_{i,j}}\big)=
h^0\big(X,\w_X(-kA-k'B)\big).
$$
Since $A$ and $B$ are general points 
of $X$, and $k+k'\leq g-2$, it follows that
\begin{equation}\label{h0F2}
h^0(\wh F_{i,j},\cal L_{i,j}|_{\wh F_{i,j}})=g-1-k-k'.
\end{equation}

Consider now the natural short exact sequence,
\begin{equation}\label{lij}
0\to\cal L_{i,j}\to\cal L_{i,j-1}\to
\cal L_{i,j-1}|_{F_{i,j}}\to 0.
\end{equation}
Restricting it to $C$, we obtain the short exact sequence
$$
0\to\cal L_{i,j}|_{\wh F_{i,j}}\to\cal L_{i,j-1}|_C\to
\cal L_{i,j-1}|_{F_{i,j}}\to 0.
$$
So, putting together \eqref{h0F} and \eqref{h0F2}, 
we get
$$
h^0(C,\cal L_{i,j-1}|_C)\leq 
h^0(F_{i,j},\cal L_{i,j-1}|_{F_{i,j}})+
h^0(\wh F_{i,j},\cal L_{i,j}|_{\wh F_{i,j}})=g.
$$
By semicontinuity, equality 
holds above. Also, each of the 
restriction maps in the composition below is surjective:
\begin{equation}\label{surjj}
H^0(\TS,\cal L_{i,j-1})\lra H^0(C,\cal L_{i,j-1}|_C)
\lra H^0(F_{i,j},\cal L_{i,j-1}|_{F_{i,j}}).
\end{equation}
As a consequence, not only is the sequence derived from 
\eqref{lij},
\begin{equation}\label{Hlij}
0\to H^0(\TS,\cal L_{i,j})\to H^0(\TS,\cal L_{i,j-1})\to 
H^0(F_{i,j},\cal L_{i,j-1}|_{F_{i,j}})\to 0,
\end{equation}
left-exact, but also right-exact.

Now, the inclusion map
$$
H^0(\TS,\cal L_{i,j})\lra H^0(\TS,\cal L_{i,j-1})
$$
is a map of free 
$\mathbb C[[t]]$-modules of the same rank 
$g$. Because of the exactness 
in \eqref{Hlij}, and because of 
\eqref{h0F}, the determinant of this map is an element of 
$\mathbb C[[t]]$ of order $k+k'+1$. Summing $k+k'+1$ for 
$j=0,1,\dots,i(g-i)-1$, we get that
the determinant of the inclusion
$$
H^0(\TS,\cal L_i)\lra H^0(\TS,\w_p)
$$
has order $i(g-i)g/2$. Thus, using \eqref{Li}, 
and using the functoriality of the ramification scheme 
(see Subsection \ref{prel2}), we get
$$
\text{mult}_{E_i}(W(p))+\frac{i(g-i)g}{2}=gi(g-i)+
\text{mult}_{E_i}(W(p_*\cal L_i,\cal L_i)).
$$
However, since the kernel of the restriction 
$H^0(\TS,\cal L_i)\to H^0(E_i,\cal L_i|_{E_i})$ is 
$tH^0(\TS,\cal L_i)$, we have that 
$\text{mult}_{E_i}(W(p_*\cal L_i,\cal L_i))=0$. Thus
$$
\text{mult}_{E_i}(W(p))=gi(g-i)-\frac{i(g-i)g}{2}=
\frac{i(g-i)g}{2},
$$
as stated. The transversality in the statement 
is equivalent to the fact that the limit Weierstrass 
divisor is reduced.
\end{proof}

\begin{remark}\label{irredlimrmk}
Keep the setup of Theorem~\ref{irredlimrmk}, but do not 
assume that $p$ is regular. If the 
singularity types of the nodes of $C$ in $\TS$ are equal, 
the same conclusions holds, with the only difference 
that we replace the $E_i$ by the $E_i^p$ defined in 
Subsection \ref{prel3}. In fact, making these 
substitutions, the proof given above works word by word.
\end{remark}

\begin{theorem} Let $g$ be a positive integer. 
Let $X$ be a general smooth curve of genus $g-1$, 
and $A$ and $B$ general points of $X$. Let 
$C$ be the nodal curve of genus $g$ obtained from 
$X$ by identifying $A$ and $B$. Then no point of $C$ is a 
limit of special Weierstrass points on a family of smooth 
curves degenerating to $C$.
\end{theorem}

\begin{proof} The statement is true if $g=1$, because 
an elliptic curve has no Weierstrass points. Suppose 
$g>1$ now. Let $p\:\TS\to S$ be a smoothing of $C$. We 
claim that the geometric general fiber has no special 
Weierstrass point. It is enough to show that, after base 
changes, blowups and blowdowns with center in the 
special fiber, the limit Weierstrass divisor is reduced. 

So, as seen in Subsection \ref{prel3}, we may 
replace $p$ by a smoothing 
$\wt p\:\wt\TS\to S$ 
whose special fiber is the curve $C$ described in 
Theorem~\ref{irredlim}, 
whose nodes have equal singularity types in 
$\wt\TS$.
Then, 
by Theorem~\ref{irredlim}, and 
Remark~\ref{irredlimrmk} thereafter, 
the limit Weierstrass divisor of $\wt p$ is reduced.
\end{proof}

\section{The special ramification locus}\label{applic}

\begin{claim}\label{srl}({\em Special ramification loci}) 
Let $\ol M_g$ be the coarse moduli space 
of stable curves of genus 
$g\geq 4$, and $\ol M_{g,1}$ that of pointed curves. 
Let $\ol\pi:\ol M_{g,1}\to\ol M_g$ be the 
forgetful map. 

Let $M_g\subset\ol M_g$ denote the open locus 
of smooth curves and $M^0_g\subset M_g$ that of 
smooth curves without nontrivial automorphisms. Set 
$M_{g,1}:=\ol\pi^{-1}(M_g)$ and 
$M_{g,1}^0:=\ol\pi^{-1}(M_g^0)$, and 
let $\pi\: M_{g,1}\to M_g$ and 
$\pi^0\: M_{g,1}^0\to M_g^0$ be the restrictions of 
$\ol\pi$. Recall that $M^0_g$ is a fine moduli space, 
with universal family $\pi^0$. Also, 
since $g\geq 4$, the boundary of 
$M_g^0$ in $M_g$ has codimension 
at least
2; 
see \cite{HaMo}, p. 53.

Let $\ol{VSW}_g$ denote the closure of 
$VSW(\pi^0)$ in $\ol M_{g,1}$; see 
Subsection \ref{prel2}. 
Also, let $\ol{VE}_{g,j}$ denote the closure of 
$VE_j(\pi^0)$ in $\ol M_{g,1}$ for $j=-1,1$.

Notice that $\ol{VSW}_g$ and the schemes $\ol{VE}_{g,j}$ 
have pure codimension 2 in $\ol M_{g,1}$. 
Indeed, their intersections with $M^0_{g,1}$ are 
determinantal of codimension at most 2. Since 
$VE_j(\pi^0)\subseteq VSW(\pi^0)$ 
set-theoretically for each $j$, it is enough to show 
that an irreducible component $U$ of $VSW(\pi^0)$ 
cannot have codimension smaller than 2. And indeed, 
since the restriction $\pi^0\: VSW(\pi^0)\to M^0_g$ 
is finite, if 
$\text{codim}(U,M^0_{g,1})\leq 1$, 
then 
$\pi^0(U)=M^0_g$. However, this is not possible 
because a general curve has no special Weierstrass 
points; see \cite{cumestgat1}, Prop. 3.1.

We observe that $\ol{VSW}_g\cap M_{g,1}$ parameterizes 
the smooth pointed curves $(C,P)$ such that $P$ is a 
special Weierstrass point of $C$. 
Indeed, there is a smooth, 
projective map $p\:\TS\to S$ 
whose fibers are curves of genus $g$, 
and such that the horizontal maps in the 
naturally induced commutative diagram
      $$\begin{CD}
        \TS @>\phi_1>> M_{g,1}\\
	@VpVV @V\pi VV\\
	S @>\phi >> M_g
	\end{CD}$$
are finite and surjective; see \cite{HaMo}, 
Lemma 3.89, p. 142. The image of 
$VSW(p)$ under $\phi_1$ parameterizes the 
smooth pointed curves $(C,P)$ such that $P$ is a 
special Weierstrass point of $C$. We need to show 
that 
$$
\phi_1(VSW(p))=\ol{VSW}_g\cap M_{g,1}.
$$
Let $S^0:=\phi^{-1}(M_g^0)$. Since 
$\phi$ is finite and surjective, 
the boundary of $S^0$ in $S$ 
has codimension 2. Reasoning as in the last 
paragraph, we see that the irreducible components of 
$VSW(p)$ meet $\TS^0:=p^{-1}(S^0)$. Thus 
$VSW(p)$ is set-theoretically the closure of 
$VSW(p^0)$ in $\TS$, 
where $p^0:=p|_{\TS^0}\:\TS^0\to S^0$. 
Then
\begin{align*}
\phi_1(VSW(p))=&\phi_1\Big(\ol{VSW(p^0)}\Big)=
\ol{\phi_1(VSW(p^0))}\cap M_{g,1}\\
=&\ol{VSW(\pi^0)}\cap M_{g,1}=\ol{VSW}_g\cap M_{g,1}.
\end{align*}

An analogous reasoning shows that
$\ol{VE}_{g,j}\cap M_{g,1}$ parameterizes
the smooth pointed curves $(C,P)$ such that 
$P$ is a special Weierstrass point of $C$ of type 
$g+j$, for $j=-1,1$. 

Set 
$\ol E_{g,j}:=\ol\pi_*[\ol{VE}_{g,j}]$ for $j=-1,1$ 
and $\ol{SW}_g=\ol\pi_*[\ol{VSW}_g]$. 
Since 
$\ol\pi|_{\ol{VSW}_g}$ is generically finite, 
$\ol{SW}_g$ and the 
$\ol E_{g,j}$ are cycles of pure codimension 1 of 
$\ol M_g$.
\end{claim}

\begin{claim}\label{pic}({\em The Picard group}) 
Let $\ol M_g$ be the coarse moduli space 
of stable curves of genus $g\geq 3$. Since $\ol M_g$ 
has only finite quotient singularities, the group of 
codimension-1 cycle classes of $\ol M_g$ 
with rational coefficients is 
isomorphic to the Picard group with rational 
coefficients, 
$\text{\rm Pic}(\ol M_g)\ox\QQ$. This group is 
freely generated by the tautological class $\lambda$ 
and the boundary 
classes $\delta_0,\delta_1,\dots,\delta_{[g/2]}$; see 
\cite{AC}, Thm. 1, p. 154 and \cite{HaMo}, 
Prop. 3.88, p. 141.
\end{claim}

\begin{theorem}\label{formSW} Let $g\geq 4$. 
The following formula holds 
in $\text{\rm Pic}(\ol M_g)\ox\QQ$:
\begin{eqnarray*}
\ol{SW}_g&=&(3g^4+4g^3+9g^2+6g+2)\lambda
-\frac{g(g+1)(2g^2+g+3)}{6}\delta_0\\
&&-(g^3+3g^2+2g+2)\sum_{i=1}^{[g/2]}i(g-i)\delta_i.
\end{eqnarray*}
\end{theorem}

\begin{proof} The above formula was shown in 
\cite{mathscand}, Thm. 5.1, p. 44,
using the method of test curves. 
Here we will show how to obtain it directly. 

Let $p\:\TS\to S$ be a flat, projective map over a 
smooth, projective curve $S$ whose fiber $X_s$ 
over each $s\in S$ is a stable curve of genus $g$. 
Assume $p$ has finitely many singular fibers, all 
of them uninodal, and that each of its nonsingular fibers 
has no nontrivial automorphisms. Also, assume the 
general fiber of $p$ is a general curve of genus $g$, 
and in particular has 
no special Weierstrass points, and each of the singular 
fibers of $p$ is general, among singular fibers of like 
nature.

For each $s\in S$ such that 
$X_s$ is singular, let $P_s$ denote the unique 
node of $X_s$, and let $k_s$ be the singularity 
type of $P_s$ in $\TS$. 
Let $S_0$ be the set of $s\in S$ such that $X_s$ is 
singular and irreducible. 
In addition, for each $i=1,\dots,[g/2]$, let 
      $$S_i:=\{s\in S\,|\,X_s\text{ contains 
        a component of genus }i\}.$$
Up to replacing $S$ by a finite covering we may assume 
that $k_s+1$ is divisible by $g$ for each $s\in S_0$.

Let $\lambda':=c_1(p_*\w_p)$, where $\w_p$ is the 
relative dualizing sheaf of $p$, and set
$$
\delta'_i:=\sum_{s\in S_i}(k_s+1)[s]
$$
for $i=0,1,\dots,[g/2]$. 

Let $p^0$ the restriction of $p$ over its smooth locus. 
Since the general fiber of $p$ has no special Weierstrass 
points, $VSW(p^0)$ is finite. Consider the zero cycle 
$SW(p^0):=p_*[VSW(p^0)]$. 
Viewing $SW(p^0)$ as a divisor class in 
$\text{Pic}(S)\otimes\QQ$, to show the statement of the 
theorem we need only show that in 
$\text{Pic}(S)\otimes\QQ$ the class of $SW(p^0)$ 
satisfies an equation similar to that of $\ol{SW}_g$, 
but with $\lambda$ and the $\delta_i$ replaced by 
$\lambda'$ and the $\delta'_i$.

By considering blowups and blowdowns, 
we may find a map of schemes 
$\beta\:\wt{\TS}\to\TS$ such that 
\begin{enumerate}
\item $\beta$ is an isomorphism away from the points 
$P_s$ for $s\in S_0$;
\item for each $s\in S_0$, the fiber 
$\wt X_s:=(p\circ\beta)^{-1}(s)$ is the 
nodal curve that is the union of the normalization of 
$X_s$ and a chain of $g-1$ rational curves connecting the 
branches over $P_s$, and $\beta\:\wt X_s\to X_s$ is the 
map collapsing the chain to $P_s$;
\item the singularity type in $\wt{\TS}$ 
of each of the nodes of $\wt X_s$ is $(k_s+1)/g-1$ for 
each $s\in S_0$.
\end{enumerate}

Let $\wt p:=p\circ\beta$. 
For each $s\in S_0$, let $\wt k_s:=(k_s+1)/g-1$, let 
$\wt X^\nu_s\subset\wt X_s$ be the normalization of 
$X_s$, and let $E_s=(E_{s,1},\dots,E_{s,g-1})$ 
be the chain of rational components of $\wt X_s$. Also, 
for each $i\geq 1$ and each $s\in S_i$, let $Y_s$ denote 
the component of the fiber $\wt X_s$ of genus $i$ and 
$Z_s$ that of genus $g-i$. Notice that $Y_s^{\wt p}$ and 
$Z_s^{\wt p}$ are Cartier divisors of $\TS$ such that 
$Y_s^{\wt p}+Z_s^{\wt p}=(k_s+1)\wt p^*(s)$ and 
$$
Y_s^{\wt p}\cdot Z_s=Z_s^{\wt p}\cdot Y_s=1.
$$
Likewise, $(\wt X^\nu_s)^{\wt p}$ and the 
$E_{s,i}^{\wt p}$ 
are Cartier divisors of $\wt\TS$ such that
$$
(\wt X^\nu_s)^{\wt p}+\sum_{i=1}^{g-1}E_{s,i}^{\wt p}=
(\wt k_s+1)\wt p^*(s),
$$
and 
$$
(\wt X^\nu_s)^{\wt p}\cdot E_{s,i}=
E_{s,i}^{\wt p}\cdot \wt X^\nu_s=
\begin{cases}
1&\text{if $i=1$ or $i=g-1$},\\
0&\text{otherwise}
\end{cases} 
$$
for each $i=1,\dots,g-1$, while for each distinct 
$i,j=1,\dots,g-1$,
$$
E_{s,i}^{\wt p}\cdot E_{s,j}=
\begin{cases}
1&\text{if $|j-i|=1$},\\
0&\text{otherwise}.
\end{cases} 
$$

Let $\wt p:=p\circ\beta$ and consider the Weierstrass 
divisor $W_{\wt p}$. Let
\begin{align*}
W:=W_{\wt p}\,-&\sum_{i=1}^{[g/2]}\sum_{s\in S_i}
\Big(\binom{g-i+1}{2}Y_s^{\wt p}+
\binom{i+1}{2}Z_s^{\wt p}\Big)\\
-&\sum_{s\in S_0}\sum_{i=1}^{g-1}
\frac{i(g-i)g}{2}E_{s,i}^{\wt p}.
\end{align*}
We claim that $W$ is effective and intersects 
transversally each singular fiber of $\wt p$.

Indeed, the claim can be checked infinitesimally around 
each $s\in S$ for which $\wt X_s$ is singular. So, it 
is possible to treat the fibers over $S_0$ and over 
the $S_i$ for $i>0$ independently.

First, Cukierman showed in \cite{cukie}, 
Prop. 2.0.8, p. 325, that $W$ is effective on a 
neighborhood of the fiber over any $s\in S_i$, for $i>0$, 
and that $W$ intersects properly this fiber. 
(In fact, Cukierman assumed $\wt p$ regular, 
but his proof goes through in our more 
general situation.) The intersection of $W$ with the 
fiber is the limit Weierstrass divisor, which 
can be computed using Formula \eqref{dR}, as in the 
proof of Theorem \ref{irredlim}, and shown 
to be reduced. So $W$ intersects transversally each 
fiber $\wt X_s$ for each $s\in S_i$ and each 
$i=1,\dots,[g/2]$.

Finally, our Theorem \ref{irredlim}, coupled with Remark 
\ref{irredlimrmk}, shows that 
$W$ is effective on a 
neighborhood of the fiber over any $s\in S_0$, and 
intersects that fiber transversally. 
This finishes the proof of the claim.  

Since $W$ intersects transversally each singular fiber 
of $\wt p$, its branch locus over $S$ is simply 
$VSW(p^0)$. Since $VSW(p^0)$ has codimension 2 in 
$\wt\TS$, we have 
$$
[VSW(p^0)]=c_2(J^1_{\wt p}(\cal O_{\wt\TS}(W)))=
c_1(W)\big(c_1(W)+c_1(\w_{\wt p})\big).
$$
In addition, $c_1(W)$ can be computed from the definition 
of $W$, since, by Pl\"ucker formula, 
$$
c_1(W_{\wt p})=\binom{g+1}{2}c_1(\w_{\wt p})-
\wt p^*\lambda'.
$$
It is now a simple but tedious task, using the 
intersection theory of $\wt\TS$, and the formula 
(see \cite{HaMo}, Formula 3.110, p. 158)
$$
\wt p_*(c_1(\w_{\wt p})^2)=12\lambda'-\delta'_0-\delta'_1
-\cdots-\delta'_{[g/2]},
$$
to compute 
$SW(p^0)=p_*[VSW(p^0)]$ and obtain the stated formula.
\end{proof}

\section{The special ramification loci of type $g+j$}

\begin{proposition}\label{schthunion} Let $g\geq 4$. 
The following 
statements hold:
\begin{enumerate}
\item $\ol{VE}_{g,j}\subseteq\ol{VSW}_g$ for $j=-1,1$.
\item Set-theoretically, 
$\ol{VE}_{g,-1}\cup\ol{VE}_{g,1}=\ol{VSW}_g$.
\item In the cycle group of codimension-$2$ 
cycles of $\ol M_{g,1}$,
    $$[\ol{VSW}_g]=[\ol{VE}_{g,-1}]+[\ol{VE}_{g,1}].$$
\end{enumerate}
\end{proposition}

\begin{proof} 
The statements are local, and can be checked 
on a neighborhood of a point 
of $M^0_{g,1}$, parameterizing a pointed stable curve 
$(C,P)$. 

Locally, the scheme $VE_1(\pi^0)$ is given 
by all maximal minors of a (Wronskian) matrix of regular 
functions of the form
        $$M=\begin{bmatrix}
            A\cr c \cr d
	    \end{bmatrix},$$
where $A$ is a matrix with $g$ columns and $g-1$ rows, 
and $c$ and $d$ are row vectors of 
size
$g$. 
Furthermore, $VE_{-1}(\pi^0)$ is given by all 
maximal minors of the matrix $A$, and $VSW(\pi^0)$ is 
given by the determinants of the square submatrices
        $$M_1:=\begin{bmatrix}
          A \cr c
          \end{bmatrix}\quad{\rm and }\quad
          M_2:=\begin{bmatrix}
	  A \cr d
	  \end{bmatrix}.$$
Since the determinants of $M_1$ and $M_2$ 
are also maximal minors of $M$, 
it is clear that $VSW(\pi^0)\supseteq VE_1(\pi^0)$. 
On the other hand, expanding these two  
determinants by the last rows, we see that they 
belong to the ideal of maximal minors of $A$. Thus 
$VSW(\pi^0)\supseteq VE_{-1}(\pi^0)$ as well. 
Statement 1 is proved.

As for Statement 2, just observe that 
if $A$ has rank $g-1$ at $(C,P)$, and thus 
$(C,P)$ is not in $VE_{-1}(\pi^0)$, then 
the vanishing of 
the determinants giving $VSW(\pi^0)$ at $(C,P)$ says 
that the two last rows of $M$ depend linearly on the 
rows of $A$, and hence $M$ has rank $g-1$, yielding that 
$(C,P)$ is in $VE_1(\pi^0)$.

Finally, for the last statement we may resort to 
Lemma 5.3 in \cite{cumestgat1}. To apply this lemma, and 
immediately get the 
last statement, we need only check that 
$M_1$ has rank at least $g-1$ at $(C,P)$, or equivalently, that $h^0(C,\omega_C(-gP))\leq 1$. 
However, it follows from \cite{Diazduke}, Thm. 4.13, 
p. 918, and Claim 3 on p.~920, that the subset of 
$M_{g,1}$ parameterizing pointed curves 
$(C,P)$ such that $h^0(C,\omega_C(-gP))\geq 2$ has 
codimension at least 3. Since the equality we want to 
prove involves codimension-2 cycles, we 
may indeed assume that $h^0(C,\omega_C(-gP))\leq 1$.
\end{proof}

\begin{claim}\label{testfamily} 
(\emph{Computing the classes $\ol E_{g,j}$}) 
We may write 
    $$\ol E_{g,j}=a_j\lambda-b_{j,0}\delta_0
      -b_{j,1}\delta_1-\cdots-b_{j,[g/2]}\delta_{[g/2]}
      \in\text{Pic}(\ol{M_g})\ox\QQ$$
for $j=-1,1$, where the $a_j$ and the $b_{j,\ell}$ 
are rational numbers to be computed. 

The coefficients $a_{-1}$ and $a_1$ were determined 
by using Porteous formula on a general family 
of smooth curves in 
\cite{Diazexc}, Thm. 4.33, p. 21 and Thm.~A1.4, p. 59:
\begin{equation}\label{as}
a_{-1}=\frac{g^2(g-1)(3g-1)}{2}\quad\text{and}\quad
a_1=\frac{(g+1)(g+2)(3g^2+3g+2)}{2}.
\end{equation}

To compute the remaining numbers, 
we use test families. 
Our first family, $p_0\:\TS_0\to S_0$, 
is constructed by taking a 
general pencil of plane cubics passing through a 
fixed point, and adding to each member of the 
pencil a general smooth curve of genus $g-1$ 
meeting the cubic transversally at the fixed point on the 
cubic and at a fixed general point on the curve of genus 
$g-1$; see \cite{HaMo}, Ex. 3.140, p. 173. 
It follows from Theorem \ref{mainthm} and 
\cite{cumestgat1}, Prop.~3.1, that 
for a nonsingular member of the 
pencil, the resulting stable curve does not 
contain any limit of special Weierstrass points. 
Thus $\int_{S_0}\ol{SW}_g\geq 0$, with strict 
inequality if and only if there is a fiber of $p_0$ 
represented by a point in the support of $\ol{SW}_g$. 
However, a quick computation, using the formula for 
$\ol{SW}_g$ in Theorem \ref{formSW}, yields 
      \begin{equation}\label{intSW}
	\int_{S_0}\ol{SW}_g=0.
      \end{equation}
So, in particular,
      $$\int_{S_0}\ol{E}_{g,j}=0\quad
        \text{for $j=-1,1$},$$
yielding the relations 
\begin{equation}\label{abs}
a_j-12b_{j,0}+b_{j,1}=0\quad\text{for $j=-1,1$.}
\end{equation}
(These
relations were obtained directly by Diaz 
\cite{Diazexc}, Lemma~7.2, 
p. 40, 
for $j=-1$, and by 
Gatto \cite{mathscand}, p.~67, for $j=1$, 
and from them Gatto concluded 
\eqref{intSW}. Here we proceeded in the opposite way.) 
It is thus enough to compute the $b_{j,i}$ 
for $j=-1,1$ and $i\geq 1$.

For each $i=1,\dots,[g/2]$, 
let $X$ be a general smooth curve 
of genus $g-i$, let $Y$ be a general smooth curve 
of genus $i$, and $B\in Y$ a general point.
Identifying the diagonal 
$\Delta\subset X\times X$ with 
$\{B\}\times X\subset Y\times X$ in the natural way, 
we get a flat, projective map $p_i\:{\cal F}_i\to X$ 
whose fiber over each $P\in X$ is the uninodal 
stable curve union of $X$ and $Y$ with $P$ and $B$ 
identified; denote by $X\cup_PY$ this fiber. Let 
$\gamma_i\:X\to\overline{M}_g$ be the induced map.  

The crux of the method is to compute the degree of the 
pullback $\gamma_i^*\ol E_{g+j}$ for $j=-1,1$. First, we 
claim that the number of 
points 
$Q\in X\cup_PY$ for all $P\in X$ which are limits of 
special Weierstrass points is a lower bound for this 
degree. 

Indeed, since $\ol E_{g,j}=\ol\pi_*[\ol{VE}_{g,j}]$, 
the support of the cycle $\ol E_{g,j}$ parameterizes 
the curves which are limits of smooth curves having 
a special Weierstrass point. So, by 
Theorem \ref{mainthm}, a curve $X\cup_P Y$ is 
parameterized in this support only if either the 
complete linear system of sections of $\w_X((i+1)P)$ or 
that of $\w_Y((g-i+1)B)$ has special ramification points, 
other than $P$ or $B$, or $P$ is a Weierstrass point 
of $X$. However, by \cite{cumestgat1}, 
Prop. 3.1, since $B$ is general, the linear system 
$H^0(Y,\w_Y((g-i+1)B))$ has no 
special ramification points other than $B$, and 
the same is true for $H^0(X,\w_X((i+1)P))$ for 
a general $P$. Thus $\gamma_i(X)$ intersects the support 
of $\ol E_{g,j}$ in finitely many points.

To compute $\gamma_i^*\ol E_{g,j}$ 
we do as follows: since 
$\ol M_g$ has finite 
quotient singularities, there is an integer 
$n>0$ such that $n\ol E_{g,j}=[D]$ 
for some effective Cartier 
divisor $D$; then 
$\gamma_i^*\ol E_{g,j}=(1/n)[\gamma_i^{-1}D]$. Since 
$\gamma^{-1}_iD$ is finite, to compute it we need only 
work infinitesimally around each $P\in X$. So, 
let $\wh p_i\:\wh{\cal F}_i\to\wh X$ denote the base 
change of $p_i$ to 
$\wh X:=\text{Spec}(\wh{\cal O}_{X,P})$, and 
$\wh\gamma_i\:\wh X\to\ol M_g$ the corresponding map. 
Let $u\:{\cal U}\to T$ be the universal deformation of 
$X\cup_P Y$. Here, $T$ is the spectrum of a ring of power 
series, whence regular. Let $\beta\:T\to\ol M_g$ be the 
induced map. Because of the universal property of $u$, 
there is a map $\alpha\:\wh X\to T$ such that 
$\wh\gamma_i=\beta\circ\alpha$. We will first describe 
$\beta^{-1}(D)$, or the cycle $[\beta^{-1}(D)]$, which 
amounts to the same because $T$ is regular.

Let $T^0\subset T$ be the open subscheme parameterizing 
the smooth fibers of $u$ without nontrivial 
automorphisms. Set ${\cal U}^0:=u^{-1}(T^0)$ and 
denote by $u^0\:{\cal U}^0\to T^0$ the induced map. The 
map $\beta$ restricts to a map $\beta^0\:T^0\to M_g^0$. 
Since the boundary of $M_g^0$ in $M_g$ has codimension 
at least 2, to compute $[\beta^{-1}(D)]$ we need only 
describe $(\beta^0)^{-1}(D\cap M_g^0)$. But 
$\ol E_{g,j}\cap M^0_g$ is equal to 
$\pi_*^0[VE_j(\pi^0)]$, 
and is a Cartier 
divisor of $M_g^0$, because $M_g^0$ is smooth. Moreover,
since both $T$ and $M_g^0$ are regular, and 
since the diagram
$$
\begin{CD}
{\cal U}^0 @>>> M^0_{g,1}\\
@Vu^0VV @V\pi^0VV\\
T^0 @>\beta^0>> M^0_g
\end{CD}
$$
is Cartesian, and the formation 
of $VE_j(\cdot)$ commutes with base change, we have
$$
(\beta^0)^*(\ol E_{g,j}\cap M^0_g)=u^0_*[VE_j(u_0)].
$$
Thus
$$
[(\beta^0)^{-1}(D\cap M_g^0)]=nu^0_*[VE_j(u_0)],
$$
and so $[\beta^{-1}(D)]=nu_*[\ol{VE_j(u^0)}]$,
where $\ol{VE_j(u^0)}$ denotes the closure of 
$VE_j(u^0)$ in $\cal U$. It follows that
$$
\wh\gamma_i^*\ol E_{g,j}=(1/n)[\wh\gamma^{-1}_i(D)]=
(1/n)[\alpha^{-1}\beta^{-1}(D)]=
\alpha^*u_*[\ol{VE_j(u^0)}].
$$
Finally, since $u$ is universal, 
we have a Cartesian diagram,
$$
\begin{CD}
\wh{\cal F}_i @>\alpha_1>> {\cal U}\\
@V\wh p_iVV @VuVV\\
\wh X @>\alpha >> T,
\end{CD}
$$
that shows that 
\begin{equation}\label{ga}
\wh\gamma_i^*\ol E_{g,j}=
\wh p_{1*}\alpha_1^*[\ol{VE_j(u^0)}].
\end{equation}

Thus, it follows from \eqref{ga} that the multiplicity 
of $\wh\gamma_i^*\ol E_{g,j}$ (at the closed point of 
$\wh X$) is at least the number of points $Q\in X\cup_PY$ 
whose image in ${\cal U}$ under $\alpha_1$ lies in 
$\ol{VE_j(u^0)}$. Now, since the 
singularities of $\ol M_g$ are quotient, 
$\wh{\cal O}_{T,0}$ 
is a finite $\wh{\cal O}_{\ol M_g,C}$-module, where 
$C:=X\cup_P Y$, and $0$ is the closed point of $T$. Thus 
$T^0$ has codimension 2 in $T$ as well. Reasoning as in 
Subsection~\ref{srl}, we can show that, as a set, 
$\ol{VE_j(u^0)}$ is also the closure of the analogous 
subscheme defined for the subfamily of $u$ consisting 
of all smooth fibers. So the multiplicity 
of $\wh\gamma_i^*\ol E_{g,j}$ is at least the number of 
points $Q\in X\cup_PY$ that are limits of special 
Weierstrass points. Our claim is thereby proved.

For each $j=-1,1$, let 
$d_{j,i}$ be the number of points $(P,Q)\in X\times X$ 
such that $Q\neq P$ and $Q$ 
is a special ramification point of type 
$g+j$ of the complete linear system of sections 
of $\w_X((i+1)P)$. Also, let 
$d'_{j,i}$ be the number of ramification points different 
from $B$ of 
the complete linear system of sections of 
$\w_Y((g-i+1+j)B)$, and $d''$ be 
the number of Weierstrass points of $X$. For 
each $j=-1,1$, set
    $$e_{j,i}=d_{j,i}+d''d'_{j,i}.$$
By Theorem \ref{mainthm}, the number of 
points 
$Q\in X\cup_PY$, for all $P\in X$, 
that are limits of special 
Weierstrass points is $e_{j,i}$. Thus
      \begin{equation}\label{>eji}
	\int_X\gamma_i^*\ol E_{g,j}\geq e_{j,i}.
      \end{equation}

Now, since $Y$ and $B$ are general, 
by \cite{cumestgat1}, Prop. 3.1, the complete 
linear system of sections of $\w_Y((g-i+1+j)B)$ 
has no special ramification points, other than $B$ with 
weight $i$. 
Also, since $X$ is general, $X$ does not have any 
special Weierstrass points; see \cite{cumestgat1}, 
Cor. 3.3. Thus, by Pl\"ucker formula,
$$
d'_{j,i}=(g+j)(g+i+j-1)+(i-1)(g+j)(g+j-1)-i=i(g+j)^2-i,
$$
and 
$$
d''=(g-i-1)(g-i)(g-i+1).
$$
In addition, by \cite{cumestgat1}, Thm. 5.6,
$$
d_{j,i}=(g-i)(g-i-1)\Big((g+j)^2(i+1)^2-(g-i+j)^2\Big).
$$
Thus
$$
e_{j,i}=i(g-i)(g-i-1)\Big((g+j)^2(g+3)+2(g+j)-(g+1)\Big).
$$

Now, using Theorem \ref{formSW}, and that 
$\gamma_i^*\delta_j=0$ for every $j\neq i$, while 
$$
\gamma_i^*\lambda=0\quad\text{and}\quad
\int_X\gamma_i^*\delta_i=2(1-g+i),
$$ 
a simple computation yields:
      \begin{equation}\label{eji}
	\int_X\gamma_i^*\ol{SW}_g=e_{-1,i}+e_{1,i}.
      \end{equation}
Using Proposition~\ref{schthunion}, and using 
\eqref{>eji} for $j=-1,1$ and 
\eqref{eji}, we get 
$$
e_{-1,i}+e_{1,i}=\int_X\gamma_i^*\ol{SW}_g
=\int_X\gamma_i^*\ol E_{g,-1}+\int_X\gamma_i^*\ol E_{g,1}
\geq e_{-1,i}+e_{1,i}.
$$
Thus $\int_X\gamma_i^*\ol E_{g,j}=e_{j,i}$ for $j=-1,1$, 
and hence we may get $b_{j,i}$:
\begin{equation}\label{bs}
b_{j,i}=i(g-i)\Big((g+j)^2(g+3)+2(g+j)-(g+1)\Big).
\end{equation}
Finally, using \eqref{as}, \eqref{bs} and the relations 
\eqref{abs}, formulas for $b_{-1,0}$ and $b_{1,0}$ can 
be obtained.
\end{claim}

\end{document}